\documentclass[final,3p,times]{elsarticle}

\usepackage{amssymb}

\usepackage{mathtools}
\usepackage{booktabs}
\usepackage{makecell}
\usepackage{array}
\usepackage{longtable}
\usepackage{pdflscape}
\usepackage{caption}
\usepackage{subcaption}
\usepackage{tabularx, multirow}
\newcolumntype{M}[1]{>{\centering\arraybackslash}m{#1}}
\usepackage{color,soul} 

\usepackage[colorlinks= true]{hyperref}

\let\today\relax
\makeatletter
\def\ps@pprintTitle{%
    \let\@oddhead\@empty
    \let\@evenhead\@empty
    \def\@oddfoot{\footnotesize\itshape
         {Submitted preprint} \hfill\today}%
    \let\@evenfoot\@oddfoot
    }
\makeatother

\biboptions{sort&compress}
\begin{document}

\begin{frontmatter}



\title{A review on recent advances in scenario aggregation  methods  \\ for power system analysis}

\author[]{Aiusha Sangadiev}

\author[]{Alvaro Gonzalez-Castellanos}

\author[]{David Pozo\corref{cor}}
\ead{davidpozocamara@gmail.com}

\cortext[cor]{Corresponding author.}


\begin{abstract}
 Worldwide commitments to net zero greenhouse emissions have accelerated investments in renewable energy resources.
The requirements for operating and planning power systems are becoming stringent because of the need to take into account the uncertainty associated with renewable generation. 
Several modeling frameworks that consider the inherent uncertainty in the operation and planning of the power system have been extensively studied.
Stochastic optimization has been the most popular approach among these frameworks due to its intuitive representation, especially when formulated using discrete probabilistic scenarios to represent the random variables. %
Although many scenarios representing all possible uncertain operating conditions would be needed to accurate evaluate stochastic operation and planning models, the size of the scenario set impacts computational complexity, posing a significant tradeoff between uncertainty detail representation and computational tractability. 

During the last decade, a large body of research has focused on developing new scenario aggregation methods to derive reduced scenario sets that show properties similar to the original scenario set while decreasing computational burden.  
This review provides an up-to-date, comprehensive classification and analysis of the literature related to scenario aggregation methods for addressing power system optimization problems. First, we present a general framework and the aggregation methodologies. Then, the main studies related to temporal and spatial scenario aggregation are described, followed by a bibliometric analysis of the main publication sources, authors and application problems. Finally, we provide a numerical analysis and discuss 16 aggregation methods for the transmission expansion planning problem. Finally, recommendations, opportunities and conclusions are discussed.

\end{abstract}




\end{frontmatter}


\section{Introduction}
\label{sec:Introduction}

Electrical power systems are an important part of modern infrastructure that provides a steady flow of electrical energy to a variety of customers, from small residential units, where electricity is used to power home appliances, to industrial parks, where electric-powered machinery is used for large-scale manufacturing. 
The importance of electric energy supply highlights the need for power system operators to effectively and accurately model such systems for the purposes of analysis, planning, and operation optimization through the use of energy system models (ESMs) \cite{Hoffman1976}.

The first ESMs used computer-aided modeling and emerged in the 80s owing to the increase in available computational resources \cite{Lopion2018}. 
Models such as the Brookhaven Energy System Optimization Model (BESOM) \cite{Kydes1980} and its extensions MARKAL \cite{Kydes1981} and TESOM \cite{Fishbone1981} mostly focused on the evaluation of energy technologies and policies, as well as long-term technological planning. 
The adoption of the Kyoto Protocol in 1997 shifted the research focus to meeting climate goals and reducing greenhouse gas emissions \cite{IGBPTerrestrialCarbonWorkingGroup1998}, which in turn led to similar developments for the ESMs: a focus on the integration of renewable energy sources (RESs) and the increased flexibility of energy systems \cite{Pfenninger2014}. 
Because the major renewable sources, solar and wind power, are intermittent and heavily depend on weather patterns, the complexity of ESMs increased rapidly, as they needed high spatial and temporal resolution for their accurate representation \cite{Haller2012, AngelisDimakis2011, Lopion2018}. 
This increase in complexity has been somewhat offset by remarkable progress in the semiconductor industry \cite{Robison2012} and development of more effective optimization software \cite{Koch2013}, but the computational gains have been slowing down in the recent years \cite{Koch2013}.
There is a need to find alternative solutions to the growing complexity of ESMs while maintaining an acceptable level of modeling accuracy.

This issue can be addressed by reducing the size of the ESM by grouping similar objects, like network nodes or RES profiles, into a single entity. 
The downside of this approach is that it inevitably leads to a loss of information and decreases modeling accuracy. However, in most cases, the improvement of computational time is significant enough to offset accuracy losses. 
As a result, aggregation methods have been widely applied in literature on a variety of energy system problems: unit commitment (UC) \cite{Palmintier2014, Wogrin2014, Poncelet2016, Shayesteh2016, TejadaArango2018, Teichgraeber2019, Mallapragada2018}, generation expansion planning (GEP) \cite{Merrick2016, Liu2018, Mallapragada2018, Buchholz2019, Yeganefar2020, Helist2020, Domnguez2021, Gonzato2021}, transmission expansion planning (TEP) \cite{Fitiwi2015, Alvarez2017, Ploussard2017, MajidiQadikolai2018}, joint GEP and UC \cite{Palmintier2011, Scott2019, Tso2020}, joint GEP and TEP \cite{Sun2019}, long-term energy system planning \cite{Poncelet2016}, optimal power flow (OPF) \cite{Shayesteh2016}, integrated heat-and-power operation models \cite{vanderHeijde2019}, power system dimensionality reduction \cite{Rsnen2009, SanchezGarcia2014, Vankov2020}, system analysis \cite{SatreMeloy2020, Li2020, Tanoto2020} and real-time decision making \cite{Zhang2021}.

Aggregation methods can be used to reduce the dimensionality of ESMs via the spatial dimension, i.e. generation units \cite{Palmintier2011, Palmintier2014}, or network nodes \cite{SanchezGarcia2014, Shayesteh2016} or via the temporal dimension, i.e. renewable generation and demand scenarios \cite{Merrick2016, Poncelet2016, Shayesteh2016, Liu2018, MajidiQadikolai2018, Mallapragada2018, TejadaArango2018, Buchholz2019, vanderHeijde2019, Scott2019, Sun2019, Helist2020, SatreMeloy2020, Vankov2020, Tso2020, Yeganefar2020, Domnguez2021, Gonzato2021, Zhang2021}, operational states \cite{Wogrin2014, Fitiwi2015, Ploussard2017}, generation portfolios \cite{Tanoto2020}, decarbonization pathways \cite{Li2020}, electricity prices \cite{Teichgraeber2019}. Heuristic \cite{Palmintier2011, Palmintier2014, Shayesteh2016} and spectral clustering \cite{SanchezGarcia2014, Vankov2020} algorithms are mostly used for spatial dimensionality reduction, while there is a greater range of choices in reducing temporal data: heuristic algorithms,  \cite{Poncelet2016, Merrick2016, Mallapragada2018, Buchholz2019}, clustering and scenario reduction algorithms \cite{Wogrin2014, Fitiwi2015, Shayesteh2016, Ploussard2017, Liu2018, MajidiQadikolai2018, Mallapragada2018, TejadaArango2018, Buchholz2019, Scott2019, Sun2019, Helist2020, SatreMeloy2020, Vankov2020, Domnguez2021, Gonzato2021}, self-organizing maps \cite{Tanoto2020, Yeganefar2020}, optimization-based scenario selection \cite{Buchholz2019, vanderHeijde2019, Gonzato2021}, adaptive/iterative methods \cite{Scott2019, Zhang2021}.

While there has been a substantial number of reviews on the aggregation of temporal information in relation to energy system problems \cite{Ehsan2019, Hoffmann2020, Hoffmann2021}, spatial aggregation has received comparatively less attention \cite{MartnezGordn2021}, and to our knowledge there has been no comprehensive survey which included both. Moreover, while the existing reviews provide useful information and guidance regarding the details and application of the methods, direct comparison between them is often not presented due to the differences in, among other things, how the ESM is formulated, computational resources, and method implementation. Therefore, this survey paper aims to provide:
\begin{itemize}
    \item a detailed description of and introduction to the most commonly used spatial and temporal aggregation methods;
    \item a comprehensive review of the literature on how these methods are applied to problems in the energy sector;
    \item an avenue for the direct comparison of methods by conducting numerical experiments on the transmission expansion planning (TEP) problem using the most popular algorithms.
\end{itemize}

This review is structured as follows. First, Section \ref{sec:Framework} covers the aggregation framework that is commonly used in the aggregation methods literature, as well as data normalization methods, similarity measures, and performance metrics. Next, Section \ref{sec:Aggregation} provides a detailed description of methods and reviews their application to power system problems. In Section \ref{sec:Bibliometrics}, we provide a bibliometric analysis of the literature reviewed. Section \ref{sec:Experiments} presents the numerical experiments for a set of aggregation methods and addressing the transmission expansion planning problem with two network models, the Garver 6-bus network \cite{Garver1970} and a high renewable penetration (HRP) 38-bus system \cite{Zhuo2020}. Finally, Section \ref{sec:Conclusion} provides concluding remarks and a practical guide for selecting aggregation methods for ESMs.

\section{Aggregation framework}
\label{sec:Framework}

The different methodologies currently employed for information aggregation in power system models are highly varied in nature and application. However, most of them follow a similar general framework, which can be summarized in the following steps \cite{Mallapragada2018, Teichgraeber2019, Hoffmann2020}:
\begin{enumerate}
    \item Normalization (\textit{optional});
    \item Aggregation;
    \item Representation of groups through a single object;
    \item Weighting and re-scaling;
    \item Evaluation of aggregation quality (\textit{optional}).
\end{enumerate}

The main steps and their characteristics for the general aggregation framework are described in Figure \ref{Fig: AggregationMethods}. 
In the following subsections, we are going to discuss the details of steps 1, 3, 4, and 5, while aggregation methods will be covered in Section \ref{sec:Aggregation}.

\begin{figure}
    \centering
    \includegraphics[width=\linewidth]{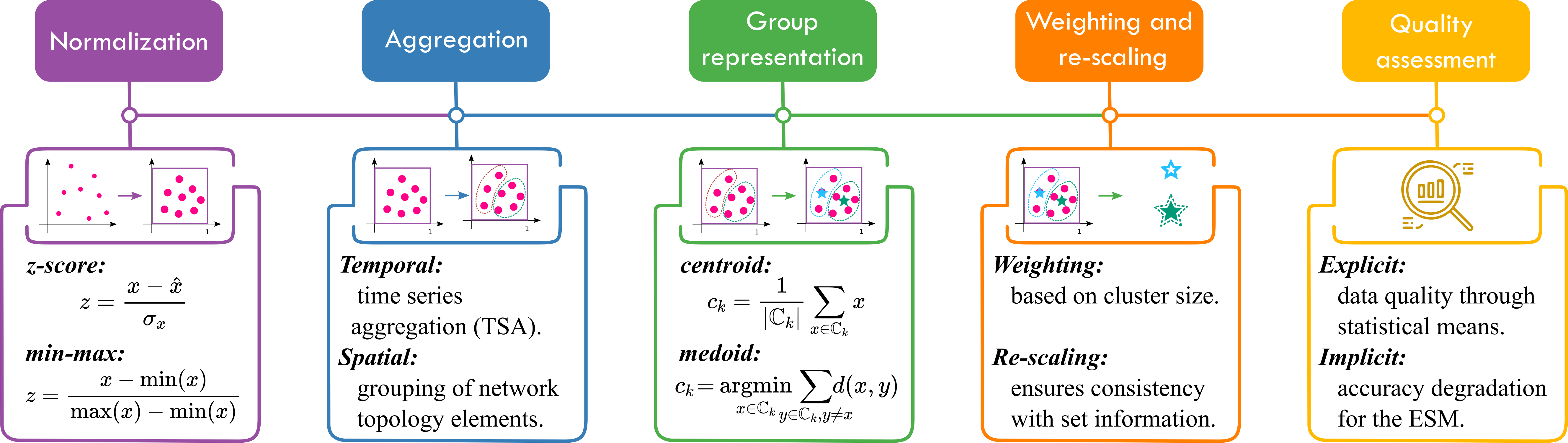}
    \caption{Sequential steps in the generalized aggregation framework.}
    \label{Fig: AggregationFramework}
\end{figure}

\subsection{Normalization}

The purpose of normalization is to scale raw data to lie within a given range, for example in $[0,1]$. This is useful if the data has different units of measurement or the values are large enough to slow down the subsequent calculations. In the literature, three normalization methods predominate: \textit{z-score} \cite{Liu2018, Teichgraeber2019, Scott2019}, \textit{min-max} \cite{Mallapragada2018, Sun2019, Tanoto2020} and division by absolute maximum value \cite{Fitiwi2015, Poncelet2016, Domnguez2021}. While the latter is self-explanatory, the other two need to be defined.
The z--score can be calculated as:
\begin{equation} \label{eq: z-score}
    z = \frac{x - \hat{x}}{\sigma_x},
\end{equation}
where $z$ is the normalized data vector, $x$ the raw data vector, $\hat{x}$ its mean value, and $\sigma_x$ its standard deviation. 

The min-max normalization is defined as:
\begin{equation} \label{eq: min-max}
    z = \frac{x - \text{min}(x)}{\text{max}(x) - \text{min}(x)}.
\end{equation}

In some papers \cite{TejadaArango2018}, normalization was used but the method was not specified or a more exotic scaling method was applied \cite{Li2020, SatreMeloy2020}. 
In principle, the choice of a specific scaling method should not affect the aggregation process \cite{Hoffmann2020}.

\subsection{Group representation}
The aggregation step groups similar objects into distinct \textit{clusters}, subdividing or partitioning the initial set of objects. Since the cluster members are similar to each other, they can be replaced by a single \textit{representative object} \cite{Teichgraeber2019, Hoffmann2020}. This step could either be integrated into the aggregation procedure, for example in the form of scenario reduction, optimization-based scenario selection, or heuristic algorithms, or applied as a separate procedure over the obtained groupings. In the latter case, usually either the cluster \textit{centroid} \cite{Rsnen2009, Fitiwi2015, Liu2018, MajidiQadikolai2018, Buchholz2019, Teichgraeber2019, Scott2019, SatreMeloy2020, Li2020, Tanoto2020, Domnguez2021} or \textit{medoid} \cite{TejadaArango2018, Teichgraeber2019, Sun2019, Gonzato2021} is chosen as the representation. In some studies, the object closest to the cluster centroid \cite{Ploussard2017, Mallapragada2018, Helist2020} is used instead, and there are also cases where the exact representation method is not specified \cite{Vankov2020, Tso2020}.

The centroid is the mean of all objects within the cluster and can be defined through the following formula:
\begin{equation} \label{eq: centroid}
    c_k = \frac{1}{|\mathbb{C}_k|} \sum\limits_{x \in \mathbb{C}_k} x,
\end{equation}
where $c_k$ is the representation of cluster $k$, $\mathbb{C}_k$ the set of objects in cluster $k$, $x$ the object within the cluster. Sometimes, the centroid cannot be calculated or use of an existing cluster object is preferred. In this case, it is suitable to represent the cluster through its medoid--the object that is the least dissimilar to all other objects within the cluster. 
The mathematical formulation for obtaining the medoid is presented below:
\begin{equation} \label{eq: medoid}
    c_k = \underset{x \in \mathbb{C}_k}{\text{argmin}} \sum\limits_{y \in \mathbb{C}_k, y \neq x} d(x,y),
\end{equation}
where $d(x,y)$ is a similarity measure between objects $x$ and $y$.

The authors in \cite{Teichgraeber2019} note that, in applying time-series aggregation to uncertainty scenarios, using centroid representations for the reduced scenario set leads to optimization problems where a lower bound is provided for the optimal solution and increasing the reduced set size shifts this bound closer to the optimal solution. 
In addition, the authors found that centroids lead to a more predictable representation of the operational parts of the optimization problem, while medoids better capture the variability in investment problems. 
Therefore, unlike the normalization methods, the choice of group representation method depends on the problem structure.

\subsection{Weighting and re-scaling}

After selecting the representative scenario from the cluster, its values need to be weighted according to the size of the cluster. If the cluster representation is a centroid, then the standard weighting will be consistent with the overall values from the full set. However, if the representation is chosen as an object from the cluster, then consistency is not guaranteed and additional re-scaling is necessary \cite{Hoffmann2020, Hoffmann2021}. 

One of the most popular approaches to re-scaling is the one first presented by Nahmmacher et al. \cite{Nahmmacher2016}, where representative scenarios were scaled by the overall annual average and its variations \cite{Mallapragada2018, Teichgraeber2019, Helist2020}. Scott et al. \cite{Scott2019} re-scale by solving a quadratic programming optimization problem, such that weighted scenarios would meet target annual net demand and peak net demand in all years. Finally, some studies use non-centroid representations but omit re-scaling \cite{Ploussard2017, TejadaArango2018, Gonzato2021}.

\subsection{Performance metrics}

Metrics for assessing aggregation quality can be divided into two groups - \textit{explicit} and \textit{implicit}. Explicit metrics directly evaluate the quality of data division through statistical means; the most commonly used metrics include the \textit{Silhouette score} \cite{Vankov2020, Li2020, SatreMeloy2020}, the \textit{Davies-Bouldin index} \cite{Li2020, SatreMeloy2020, Tanoto2020} and the \textit{Calinski-Harabasz score} \cite{Li2020, SatreMeloy2020, Tanoto2020}.
While these metrics explicitly assess the quality of aggregation, they offer limited insight onto how the partitioned and reduced data will affect the accuracy of the ESM, which is often the main concern.

Keeping the above in mind, many studies instead use implicit metrics,  - they do not directly evaluate the quality of aggregation but try to capture the degradation in accuracy for the ESM with a reduced scenario set. These metrics are often problem-specific and thus vary from study to study, but a more general metric that is found in a variety of papers is the difference in objective function values for problems solved with the full data set versus the reduced data set\cite{Palmintier2011, Palmintier2014, Wogrin2014, Fitiwi2015, Shayesteh2016, Ploussard2017, MajidiQadikolai2018, TejadaArango2018, Scott2019, Sun2019, Teichgraeber2019, Helist2020, Tso2020, Yeganefar2020, Domnguez2021, Gonzato2021}.

For two-stage stochastic optimization problems, there is an additional distinction between \textit{in-sample} and \textit{out-of-sample} objective values. 
For the former, both the first and second stages are computed with reduced data, while the latter uses first-stage decisions obtained with the reduced set and evaluates the second stage with the complete data set to assess optimality.

The main downside of implicit metrics lies in the need to compute values with full data set for comparison, which can be often computationally expensive. However, these metrics give a direct insight into the effect of aggregation on the ESM.

\section{Aggregation methods}
\label{sec:Aggregation}

Aggregation methods can be divided into two groups by the nature of data they are applied to, i.e., spatial and temporal. 
Within these groups, the methods can be further subdivided into categories based on their working principle, e.g., clustering methods or scenario reduction algorithms.
Table \ref{table:Methods} provides an overview and classification of aggregation methods, the energy problem applications, the performance metrics used to assess them, and the aggregation objects.

\begin{figure}
    \centering
    \includegraphics[width=0.95\linewidth]{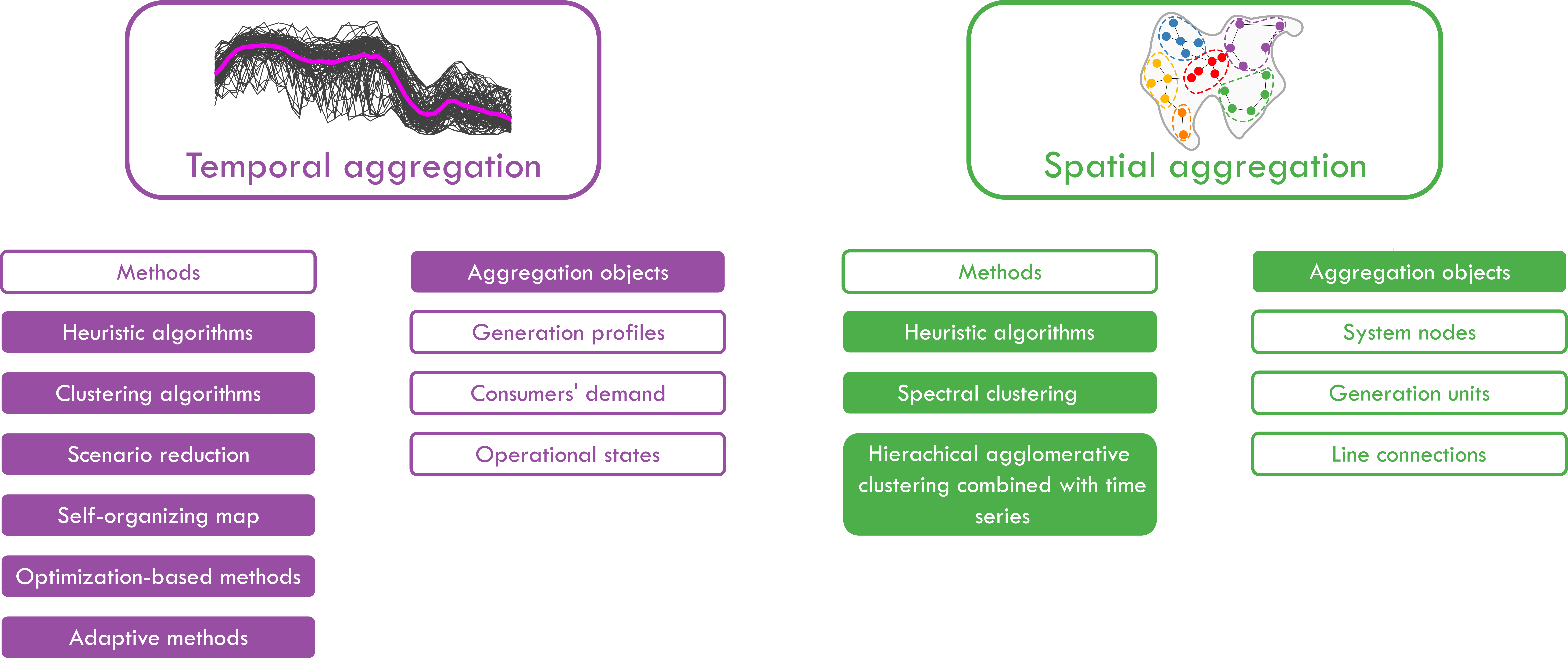}
    \caption{Summary of the main aggregation methods in energy system models.}
    \label{Fig: AggregationMethods}
\end{figure}

A summary of the main aggregation methods is presented in Figure \ref{Fig: AggregationMethods}.

In the following subsections, we describe the commonly used methods and provide an overview of the respective works that applied them.

\subsection{Temporal aggregation}

Time-series aggregation (TSA) methods have been proposed to reduce the temporal dimensionality of ESM. Time series are characterized by high dimensionality \cite{Lin2004}, noise \cite{Lin2003_1}, and often times differing lengths. The most common sources of the time-series data in ESMs are short- and long-term uncertainties linked to renewable energy generation and demand. Short-term uncertainties include conditional variability of demand or renewable generation, while an example of long-term uncertainty would be the expected growth of demand or renewable energy integration \cite{Velloso2020}.

Aggregation requires that the measure of similarity between the time series be defined first.
There are three main approaches to defining the similarity of time series \cite{Aghabozorgi2015, Vankov2020}. Two time series are similar in:
\begin{itemize}
    \item \textit{time} if their values at each time step are similar;
    \item \textit{shape} if they have similar characteristics that are not overlapping in time - e.g. they are time-shifted or scaled;
    \item \textit{structure} if they have similar laws of change - periodicity and trends.
\end{itemize}

There is no definitive guide for choosing either approach, since it largely depends on the data and task at hand \cite{Aghabozorgi2015}.

\begin{table}[htb]
    \centering
    
    \caption{Distance metrics for TSA}
    \label{table: Distance metrics}
    
    \resizebox{\textwidth}{!}{
    \begin{tabular}{M{0.2\linewidth}M{0.15\linewidth}M{0.45\linewidth}M{0.2\linewidth}}
        \toprule
        \textbf{Distance metric} & \textbf{Similarity} & \textbf{Mathematical formula} & \textbf{Works} \\
        \midrule
        Minkowski & Time & $d(x, y) = \left( \sum\limits_{i = 1}^n (x_i - y_i)^p \right)^\frac{1}{p}$ & \cite{Mallapragada2018, Vankov2020, Liu2018, Rsnen2009, Helist2020, Sun2019, Domnguez2021, Li2020, Tso2020, Scott2019, SatreMeloy2020, Tanoto2020, Teichgraeber2019, vanderHeijde2019, Gonzato2021, MajidiQadikolai2018, Buchholz2019, Ploussard2017, Alvarez2017, Wogrin2014, TejadaArango2018} \\
        \midrule
        Dynamic time warping \cite{Berndt1994} & Shape & $d(x, y) = \min \sqrt{\sum\limits_{i = 1}^n w_i}$, where $W = \{ w_i \}$ - \textit{warping path} \cite{Hoffmann2020} & \cite{Teichgraeber2019, Buchholz2019, Liu2018, Vankov2020} \\
        \midrule
        Minimum jump cost \cite{Serr2012} & Structure & $d(x, y) = \sum\limits_{i = 1}^n c_i^{\text{min}}$, where $c_i^{\text{min}}$ - jump distance to the nearest object & \cite{Vankov2020} \\
        \midrule
        Shape-based \cite{Paparrizos2015} & Shape & $d(x, y) = 1 - \max \text{NCC}_c$, where $\text{NCC}_c$ - normalized cross correlation & \cite{Teichgraeber2019}\\
        \bottomrule
    \end{tabular}}
\end{table}

Table \ref{table: Distance metrics} classifies the mathematical formulae for the most commonly used similarity measures or \textit{distance metrics}, where $x$ and $y$ are two $n$-dimensional vectors. The most popular form of Minkowski distance metric is under $p = 2$, i.e., the Euclidean distance, and it is one of the more commonly used similarity measures overall. 
Other Minkowski metrics used are $p = 1$, the Manhattan distance, and $p = \infty$, the Chebyshev distance. 
Dynamic time warping (DTW) seeks to minimize the total warping path, which is a path of minimal deviations in the matrix of cross-deviations between $x$ and $y$ \cite{Paparrizos2015}. 
Shape-based distance (SBD) employs a  similar approach to DTW but maximizes the cross-correlation between time periods \cite{Paparrizos2015}. Finally, minimum jump cost (MJC) strives to find a minimum path that ``jumps'' from one time series to another on each step \cite{Serr2012, Vankov2020}.

TSA methods can be divided into six groups as a function of their working principles: heuristics, clustering algorithms, scenario reduction, neural networks, optimization-based methods and adaptive methods. In the following subsections, we describe and review methods from each group.

\subsubsection{Heuristic methods}

Heuristic algorithms mainly use rules or expert knowledge (EK) to find and aggregate similar time series.
The most simple examples are dummy selection \cite{Buchholz2019} and random sampling \cite{Helist2020}, where the former selects every $j$-th scenario from the set, and the latter picks them randomly. 
These methods are often used as baseline solutions for comparison with more sophisticated methods.

Another relatively simple approach, which has also been found in the literature \cite{Merrick2016, Poncelet2016, Mallapragada2018}, is selecting time slices based on seasonality, time of day, day of the week, and other patterns.
However, Poncelet et al. \cite{Poncelet2016} and Mallapragada et al. \cite{Mallapragada2018} found that time slices often offer lower temporal resolution compared to the selection of representative days, which has negative effects on the accuracy of the ESM.
This finding is consistent with the work of Merrick et al. \cite{Merrick2016}, where the authors compared 4-hour time slices against the global mean as the representative object.
The use of time slices offered a good approximation of the full data of the GEP model, as they provided higher temporal resolution than did simply using one averaged value.
An idea similar to time slices was used by Munoz et al. \cite{Munoz2015}, who selected and incorporated elements from the original time series into the aggregated series in a manner that preserved statistical features of the original time series, i.e., mean and standard deviation.
Buchholz et al. \cite{Buchholz2019} found that this relatively simple approach can achieve very good results and perform better than using more complex methods on the GEP model.

Residual load duration curve (RLDC) selection \cite{Koduvere2018} is a more complex heuristic selection method.
The original residual load curve (RLC) is divided into $k$ parts, and for each part the RLDC is constructed and every 13th hour, along with the first and last hours, are selected. 
The idea is to capture both the minimum and maximum load levels, as well as all different levels of yearly RLDC. 
The selected hours are then reconstructed chronologically using the chronology of one of the elements.
This method performed well when analyzing the energy system in the large-scale Balmorel model \cite{Koduvere2018} but performed very poorly in the GEP model \cite{Buchholz2019}, which the authors concluded that was due to having purely selected hours rather than time slices or representative days.

Finally, Buchholz et al. \cite{Buchholz2019} suggested dynamic blocking (DB), which builds upon the work of Ludig et al. \cite{Ludig2011}, where equal-sized time blocks were used to group the time series and their mean was used as a representation.
Instead of using static and equal-sized blocks, the authors dynamically constructed these blocks by aggregating similar elements, where similarity was determined by simple rules.
For hourly elements, their similarity to elements in the blocks was determined by checking whether their inclusion would lead to the difference between the minimum and maximum elements in a block being higher than a set threshold.
For time slices, similarity was measured by comparing the variance within that time slice to the variance of the block; if the variances were similar, then a time slice was eligible for inclusion.
This approach performed worse than the statistical representation method of Munoz et al. \cite{Munoz2015} but better than the RLDC.

\subsubsection{Clustering algorithms}

Time-series clustering is an unsupervised process of partitioning of time-series data into groups based on the selected similarity measure \cite{Rokach2009, Aghabozorgi2015}. The most commonly used methods in the literature are \textit{partitional} and deterministic \textit{hierarchical} clustering \cite{Teichgraeber2019, Hoffmann2020}.

Partitional clustering methods solve the optimization problem of minimizing the within-cluster sum of squared distances \cite{Teichgraeber2019, Hoffmann2020}:
\begin{equation} \label{eq: partitional clustering}
    \text{min} \sum\limits_{k = 1}^K \sum\limits_{x_p \in \mathbb{C}_k} d(x_p, c_k)^2,
\end{equation}
where $K$ is the number of clusters and $d(x_p, c_k)$ the distance between object $x_p$ and cluster center $c_k$. The centers are updated using the following procedure:
\begin{equation} \label{eq: center update}
    c_k = \underset{z}{\text{min}} \sum\limits_{x_p \in \mathbb{C}_k} d(x_p, z)^2,
\end{equation}
where $z$ is the cluster center that minimizes the within-cluster distance. This problem is NP-hard and the clustering algorithms using this approach often employ greedy heuristic algorithms to efficiently solve it. 

\textit{K-means} is the most widely used partitional clustering algorithm \cite{Mallapragada2018, Liu2018, Helist2020, Rsnen2009, Li2020, SatreMeloy2020, Tanoto2020, Teichgraeber2019, MajidiQadikolai2018, Buchholz2019, Ploussard2017, TejadaArango2018, Wogrin2014, Fitiwi2015}, here $d(x_p, c_k)$ is the Euclidean distance, and the cluster center $c_k$ is represented by the centroid. 
K-means is usually solved using the heuristic Lloyd's algorithm \cite{Lloyd1982, MacQueen1967}, while initial cluster centers are determined randomly or through quasi-random initialization methods like k-means++ \cite{Arthur2007}. The algorithm performs best when data clusters are hyperspherical \cite{Tan2019}.

Some of the interesting papers that have used k-means include Räsänen et al. \cite{Rsnen2009} and Fitiwi et al. \cite{Fitiwi2015}, where the authors did not use the raw time series and instead engineered features to represent them. 
In the first study, a total of 7 features were extracted from the time series and were later used for clustering by k-means: the standard deviation, mean, skewness, kurtosis, chaos, energy, and periodicity. The paper by Fitiwi et al. \cite{Fitiwi2015} went further in including additional features such as the hourly generation-demand scenarios and the selected lines with more relevant congestions and losses. 
These moments were then used as features for the k-means algorithm, with medoid representation for the clusters.

Another popular alternative is \textit{k-medoids} (used in \cite{TejadaArango2018, Teichgraeber2019, Scott2019}), which uses medoids as cluster centers $c_k$ and, unlike k-means, it allows the use of arbitrary distance metrics. 
The standard algorithm for k-medoids is called partitioning around medoids (PAM) \cite{Kaufman1990}. It is a greedy heuristic algorithm, and improved alternatives with a faster runtime, like FastPAM \cite{Schubert2019} and FasterPAM \cite{Schubert2021}, have also been developed. It has been observed that k-medoids better capture the intra-period variance of renewable generation and demand scenarios, while k-means tend to underestimate the extreme periods due to averaging \cite{Teichgraeber2019, Hoffmann2020}.

Relatively less-employed algorithms are \textit{k-shape} \cite{Paparrizos2015} and \textit{dynamic time warping barycentric averaging} (DBA) \cite{Petitjean2011}, which were used by Teichgraeber et al. \cite{Teichgraeber2019}. 
K-shape uses SBD, and the cluster update in Eq. \eqref{eq: center update} takes the form of maximization of the Raleigh Quotient, which can be solved analytically, making the update computationally efficient. 
DBA employs DTW as a distance measure and uses barycentric averaging to find centroids for the DTW metric. 
Teichgraeber et al. \cite{Teichgraeber2019} compared them with more conventional methods like k-means and found that k-shape significantly outperformed more popular methods on problems that heavily depend on capturing intra-day variability; DBA similarly captured the variability well but performed worse than k-shape.

Hierarchical clustering algorithms seek to build a hierarchy of clusters either by agglomerative (``bottom-up'') or divisive (``top-down'') approaches \cite{Maimon2005}. For the former, each object starts in a distinct cluster and clusters are merged as they go up the hierarchy, while in the latter, all objects start in a large cluster and splits are performed as they go down the hierarchy. 
Hierarchical agglomerative clustering (HAC) is the more commonly used of the two and has been employed for TSA in many energy system problems \cite{Alvarez2017, Liu2018, Buchholz2019, Sun2019, Teichgraeber2019, Li2020, Tso2020, SatreMeloy2020, Domnguez2021, Gonzato2021}.

\begin{table}[htb]
    \centering
    
    \caption{Linkage criteria for hierarchical agglomerative clustering (HAC)}
    \label{table: Linkage criteria}
    
    \begin{tabular}{M{0.2\linewidth}M{0.55\linewidth}M{0.2\linewidth}}
        
        \toprule
        \textbf{Linkage criterion} & \textbf{Mathematical formula} & \textbf{Works} \\
        \midrule
        Ward's \cite{Ward1963} & $D(\mathbb{C}_X,\mathbb{C}_Y) = \frac{|\mathbb{C}_X| |\mathbb{C}_Y|}{|\mathbb{C}_X| + |\mathbb{C}_Y|} ||c_X - c_Y||^2$ & \cite{Teichgraeber2019, Sun2019, Tso2020, SatreMeloy2020, Domnguez2021, Gonzato2021} \\
        \midrule
        Minmax \cite{Bien2011} & $D(\mathbb{C}_X,\mathbb{C}_Y) = \underset{x \in \mathbb{C}_X \cup \mathbb{C}_Y}{\min} \left[ \underset{y \in \mathbb{C}_X \cup \mathbb{C}_Y}{\max} d(x,y) \right]$ & \cite{Buchholz2019, Liu2018} \\
        \midrule
        Complete \cite{Sorensen1948} & $D(\mathbb{C}_X,\mathbb{C}_Y) = \underset{x \in X\mathbb{C}_, y \in \mathbb{C}_Y}{\max} d(x,y)$ & \cite{SatreMeloy2020, Li2020} \\
        \midrule
        Single \cite{McQuitty1957} & $D(\mathbb{C}_X,\mathbb{C}_Y) = \underset{x \in \mathbb{C}_X, y \in \mathbb{C}_Y}{\min} d(x,y)$ & \cite{SatreMeloy2020} \\
        \midrule
        Average \cite{Sokal1958} & $D(\mathbb{C}_X,\mathbb{C}_Y) = \frac{1}{|\mathbb{C}_X| |\mathbb{C}_Y|} \sum\limits_{x \in \mathbb{C}_X} \sum\limits_{y \in \mathbb{C}_Y} d(x,y)$ & \cite{Alvarez2017, SatreMeloy2020} \\
        \midrule
        Centroid \cite{Sokal1958} & $D(\mathbb{C}_X,\mathbb{C}_Y) = ||c_X - c_Y||^2 $ & \cite{SatreMeloy2020}\\ 
        \bottomrule
    \end{tabular}
\end{table}

Whether a pair of clusters needs to be merged or not is determined via \textit{linkage criteria}, which defines the similarity between two clusters.
Table \ref{table: Linkage criteria} presents the linkage criteria commonly found in the literature, where $D(\mathbb{C}_X,\mathbb{C}_Y)$ is the distance between clusters $\mathbb{C}_X$ and $\mathbb{C}_Y$.
Aside from cluster distance, the HAC can also be used with different distance measures for the objects themselves, with the exception of Ward's criterion, which seeks to minimize the within-cluster variance in a manner that is similar to k-means method and uses the Euclidean distance.

Generally, HAC performs well when the underlying data consists of many disjoint clusters because when two dissimilar groups merge, can have a drastic increase in within-cluster variance for the new cluster.
Mallapragada et al. \cite{Mallapragada2018} found that historical load and generation profiles rarely fit that description because the data is often not easily separable - while there are days with a single peak load and days with multiple peaks, there are also many days that have the characteristics of both profiles, which leads to the absence of significant ``jumps'' in intra-cluster variance and to the creation of a large cluster which contains the majority of objects with many small clusters around it, i.e., ``rich get richer'' behavior. 
The upsides of the method are its deterministic nature and no assumption of spherical clusters for the data \cite{Tan2019}.

HAC has also been modified in the literature specifically for the purpose of performing TSA over demand-generation scenarios. 
Pineda et al. \cite{Pineda2018} introduced the chronological time-period clustering (CTPC) algorithm, which modifies the linkage criterion by introducing adjacency, i.e., two clusters $\mathbb{C}_X$ and $\mathbb{C}_Y$ are adjacent if any element $x \in \mathbb{C}_X$ contains an hour that is consecutive or antecedent to an hour of any element $y \in \mathbb{C}_Y$. 
For adjacent clusters, the linkage criterion stays the same, while for non-adjacent clusters it becomes equal to infinity, meaning that the clusters will not be merged.
Gonzato et al. \cite{Gonzato2021} found that for GEP with energy storage, data clustered by CTPC tends to exhibit a bias toward short-term storage, while elevating long-term storage more.
Domínguez et al. \cite{Domnguez2021} expanded on the CTPC to include multiple chronological periods. First, the CTPC is run for a longer time period (for example, one day) to obtain the representative periods, then a second run of the CTPC further refines the representative periods by taking a smaller seasonality window. 
The authors found that the inclusion of multiple time horizons generally helped the algorithm better capture the generation-demand variability, while also reproducing the larger seasonal (daily, weekly) patterns.

There are also examples of combining partitional and hierarchical algorithms to enhance clustering performance. 
For instance, \cite{Liu2018} uses the kMHC procedure: k-means are used to produce clusters for the full scenario set, while HAC is applied to further refine the clusters.
The authors applied kMHC to the capacity expansion problem and found that this method yields investment decisions that are comparable to the ones obtained with the full set, outperforming the pure k-means algorithm.
Buchholz et al. \cite{Buchholz2019} also tested the kMHC and its variation using fuzzy clustering instead of k-means on another capacity expansion model,  finding that kMHC and its fuzzy variation (level-correlation (LC) clustering) similarly tend to produce better results than raw k-means, with LC yielding the best results relative to heuristic statistical representation.

Among the less popular clustering algorithms, Li et al. \cite{Li2020} explored density-based DBSCAN, Gaussian mixture model (GMM), and spectral clustering, in addition to k-means and HAC. 
DBSCAN \cite{Ester1996} treats clusters as areas of high density separated by areas of low density, where high-density areas are populated by \textit{core samples}, meaning that clusters, in that case, are core samples and a set of samples located close to them. 
Probabilistic GMM \cite{Murphy2012} in based on the expectation-maximization (EM) algorithm seeks to fit a mixture of multivariate Gaussian distributions to the distribution of samples and maximize the log-likelihood function.
This allows the GMM to have greater flexibility than k-means, which also makes an assumption about the Gaussian nature of the data's distribution but cannot re-shape the probability distributions to fit any kind of elliptical clusters \cite{Bishop2006}. 
Spectral clustering will be discussed in more detail in Section \ref{sec: Spatial aggregation}, as it has found more application in spatial aggregation than in TSA.
Helistö et al. \cite{Helist2020} introduced a generalized clustering algorithm called regular decomposition (RD), whose main purpose is to segment the data into an optimal amount of classes that would maximally reveal the redundancies in the data. 
The optimal partition could be found with the minimum description length (MDL) framework. 
The authors applied the method to the problem of generation expansion planning and found that the RD is an efficient, scalable method for selecting representative slices.

\subsubsection{Scenario reduction}

Scenario reduction methods were specifically developed to reduce the dimensionality of stochastic programming (SP) models. 
Kaut and Wallace \cite{Kaut2007} provide an overview of the most commonly used methods, while here we focus on the scenario reduction algorithms that were developed to specifically tackle the two-stage stochastic programming problems \cite{Heitsch2003, Dupacova2003, Morales2009, Bruninx2016} that are often used as state-of-the-art models for energy system planning problems like transmission expansion planning (TEP) or generation/capacity expansion models.

Dupačová et al. \cite{Dupacova2003} showed that the optimal solution of a simpler problem with a reduced scenario set will be close to the optimal solution of the original problem if the \textit{probability distance} between the two scenario sets is sufficiently small. The most commonly used probability distance metric is the \textit{Kantorovich distance}, which can be obtained by solving the \textit{Monge-Kantorovich mass transportation problem} \cite{Rachev2002}:
\begin{align} \label{eq: Kantorovich}
    \begin{split}
        D_K(Q, Q') = \min \Biggl\{ \sum\limits_{\omega, \omega'} c(\omega, \omega') \eta(\omega, \omega') : & \sum\limits_\omega \eta(\omega, \omega') = \pi_{\omega'},\\ & \sum\limits_{\omega'} \eta(\omega, \omega') = \pi_{\omega} \Biggr\},
    \end{split}
\end{align}
where $D_K(Q, Q')$ is the Kantorovich distance between probability distributions $Q$ and $Q'$, $\pi_{\omega}$ and $\pi_{\omega'}$ the scenario probabilities in the respective scenario sets $\Omega$ and $\Omega'$, $\eta(\omega, \omega')$ the joint probability function defined over $Q \times Q'$, and $c(\omega, \omega')$ the cost function. 
Eq. \eqref{eq: Kantorovich} can be rewritten in the following form for two-stage stochastic problems \cite{Dupacova2003}:
\begin{equation} \label{eq: Kantorovich_EQ}
    D_K(Q, Q') = \sum\limits_{\omega \in \Omega \setminus \Omega'} \pi_{\omega} \underset{\omega'}{\min} \text{ } c(\omega, \omega').
\end{equation}

Since the problem defined in Eq. \eqref{eq: Kantorovich_EQ} is NP-hard, it is usually solved using greedy heuristic algorithms such as the \textit{forward selection algorithm} (FSA) \cite{Heitsch2003, Dupacova2003}. The FSA is an iterative algorithm consisting of two main parts:
\begin{enumerate}
    \item Greedily select the representative scenarios from the set:
    $$ \underset{r \not\in \Omega_R}{\min} \sum\limits_{\omega \not\in \Omega_R \cup \{ r \}} \pi_{\omega} \underset{\omega' \in \Omega_R \cup \{ r \}}{\min} c(\omega, \omega'), $$
    where $\Omega_R$ is the reduced scenario set, and $r$ the representative scenario.
    \item Redistribute the probabilities for the representative scenarios:
    $$ q_{\omega} = \pi_{\omega} + \sum\limits_{\omega' \in \Omega \setminus \Omega_R} \pi_{\omega'}. $$
\end{enumerate}

The choice of cost function $c(\omega, \omega')$ is important as it directly affects the process of scenario selection. Heitsch et al. \cite{Heitsch2003} proposed using a \textit{Wasserstein metric} of order $p$ as a cost function $c(\omega, \omega')$. Dupačová et al. \cite{Dupacova2003} took $p = 1$ with
\begin{equation} \label{eq: Dupacova}
    c(\omega, \omega') = |h^\omega - h^{\omega'}|,
\end{equation}
where $h^\omega$ is the realization of the stochastic vector $h$ under scenario $\omega$, and obtained a variation of the FSA with the norm of the difference between pairs of random vectors as the cost function, an approach that guarantees that the representative scenarios are close to the center of mass of the original scenario set \cite{GonzalezCastellanos2021}.
Shayesteh et al. \cite{Shayesteh2016} investigated the effect of scenario reduction on ESMs and applied the method to produce a reduced scenario set for OPF and stochastic UC problems, finding that the approach significantly decreases the computational time but has a relatively large effect on the accuracy depending on the problem.

Morales et al. \cite{Morales2009} adapted the method for scenario reduction in trading in markets of futures. 
Their cost function was defined as:
\begin{equation} \label{eq: Morales}
    c(\omega, \omega') = |z_\omega^{\text{DP}} - z_{\omega'}^{\text{DP}}|,
\end{equation}
where $z_\omega^{\text{DP}}$ is the objective function value of the optimization problem when the second stage is represented only by scenario $\omega$ and the first stage is fixed at the deterministic expected-value problem (DP) solution.

Finally, the work of Bruninx et al. \cite{Bruninx2016} presented a risk-neutral version of the method for the stochastic UC problem, where the FSA cost function is defined as:
\begin{equation} \label{eq: Bruninx}
    c(\omega, \omega') = |z_\omega^{\text{SS}} - z_{\omega'}^{\text{SS}}|,
\end{equation}
where $z_\omega^{\text{SS}}$ is the value of the objective function of the optimization problem when the second stage is represented only by scenario $\omega$.

\subsubsection{Self-organizing map}

The self-organizing map (SOM), introduced by Kohonen \cite{Kohonen1982}, is an unsupervised machine learning method that performs an ordered mapping of the input data into a lower-dimensional space. Essentially, the SOM is an artificial neural network (ANN) that is trained through a competitive learning framework, i.e., the ANN nodes compete with each other for the right to ``respond'' to the input data \cite{Rumelhart1986}.

Consider a $n$-dimensional sample $x \in \mathbb{X}$, where $\mathbb{X}$ is the sample set. 
The nodes of the SOM are initialized with randomly distributed weights $W_i$, hence the Euclidean distances between the input sample $x$ and each node in the map are calculated to find the node with the smallest dissimilarity to the input data – this node $v$ is called the \textit{best matching unit} (BMU). After that, the BMU and the nodes in its neighborhood are essentially moved towards the $x$ through the following weight update:
\begin{equation} 
    W_i(t + 1) = W_i(t) + \theta(v,i,t) \alpha(t) (x - W_i(t)),
\end{equation}
where $t$ is the current iteration step, $\theta(v,i,t)$ the neighborhood function that denotes restraint due to distance from the BMU, $\alpha(t)$ the learning restraint. The process is repeated until $t$ does not reach the set iteration threshold $T$ \cite{Kohonen1982, Kohonen2001}.

Tanoto et al. \cite{Tanoto2020} used a SOM to map possible secure and affordable low-emissions technology mixes and found that it produced results similar to k-means. Yeganefar et al. \cite{Yeganefar2020} performed extraction of the extreme periods for GEP on a NREL-118 system using a SOM-based algorithm. They took two neurons from the SOM grid, grouped the instances of the selected neurons and constructed a net load duration curve (NLDC). The extreme values of this NLDC are then determined by finding values greater than the upper whisker or lower than the lower whisker of the original NLDC. DTW distance is computed between extreme values of the created NLDC and the original NLDC, and this process is repeated for each pair of neurons within the SOM to determine the pair with the lowest DTW distance.

\subsubsection{Optimization-based methods}

This class of methods formulates the problem of TSA as an optimization problem. In most cases, mixed integer linear programming (MILP) models \cite{Poncelet2017, Buchholz2019, vanderHeijde2019, Gonzato2021} are used, but there are also examples of integer programming (IP) models \cite{MajidiQadikolai2018}, mixed integer quadratic programming (MIQP) models \cite{Gonzato2021} and non-linear programming (NLP) models \cite{Zhang2021}.

Poncelet et al. \cite{Poncelet2017} introduced the following aggregation principle: first binarize the normalized duration curves (NDCs), where each bin represents the number of hours during the year surpassing a certain level of a specific attribute, then perform the same procedure for each candidate day.
After that, formulate a MILP model to minimize the sum of absolute differences between the share of time when the original NDC was out of the respective bin's range and the hours when the approximated NDC surpassed the respective bin's borders. The authors used this method to select representative days in a GEP model.
Van der Heijdge et al. \cite{vanderHeijde2019} noted that the main disadvantage of that algorithm is the loss of chronology and attempted to fix it by introducing an additional MIQP model that restores the chronology.
The resulting approach was tested in selecting representative days for the ESM of a district heating model with seasonal storage.
The authors found that the method was able to replicate the behavior of full-year optimization with good accuracy and computational time that was 11-36 times faster.

Buchholz et al. \cite{Buchholz2019} investigated three optimization-based methods: exhaustive search (ES), based on \cite{Sisternes2013}, optimized RLDC approximation (OA), \cite{Poncelet2017}, and optimized criteria search (OS), which was introduced in the work itself. 
The ES method essentially loops through every combination of the representative time periods to find $k$ periods that represent the annual RLDC, and due to the combinatorial nature of the problem, the computational time increases very quickly. The authors used weeks as time periods because even daily selection was too computationally expensive. The OS method builds on the idea of Munoz et al. \cite{Munoz2015} and formulates the problem as a MILP model: the elements are picked in such a way that they minimize the correlation between selected elements while trying to achieve variance, average levels and ramping  as close as possible to the original time series. The three methods together, along with some other heuristics and clustering methods, were applied to scenario selection in the GEP model. The algorithms performed roughly equivalently, with neither standing out too much.

Gonzato et al. \cite{Gonzato2021} also investigated three optimization-based methods in application to GEP with long-term energy storage: Poncelet et al.'s method \cite{Poncelet2017}, Buchholz et al.'s OS method \cite{Buchholz2019}, and the optimal representative days orderer (ORDO) method, which is essentially a hierarchical clustering problem formulated as an optimization problem \cite{Kotzur2018}.
The authors found that the optimization-based approaches performed worse than HAC and HAC-based CTPC methods, while Poncelet et al.'s method \cite{Poncelet2017} also required substantial time for computation.

Majidi-Qadikolai et al. \cite{MajidiQadikolai2018} proposed a dual approach, where a clustering algorithm was used to determine the initial clusters followed by optimization-based refinement, this approach is called scenario bundling. First, the authors applied k-means to cluster the data, then further refined the clusters by minimizing the dissimilarity between scenario bundles in an IP model.
This approach led to computational speeds that were 8 times faster and a better optimality gap in a large-scale TEP model.

Finally, Zhang et al. \cite{Zhang2021} proposed a decision-making oriented clustering (DMOC) framework in the form of a generalized optimization problem as a way to minimize the negative sum of the so-called performance metrics. 
The authors demonstrated that the framework can be adapted to different problems, namely, real-time pricing (RTP) and power consumption scheduling (PCS) problems, by changing the performance metric function, and they showed that the DMOC significantly reduces the optimality loss relative to the k-means algorithm.

\subsubsection{Adaptive methods}

Although some of the methods discussed above are iterative, they do not fall in the category of adaptive/iterative methods, which contain algorithms that can refine the selection scenario process via ``feedback'' information back the partition method from the ESM \cite{Tso2020, GonzalezCastellanos2021}.

Tso et al. \cite{Tso2020} explored the problem of reducing the generation-demand scenario set for the two-stage mixed-integer stochastic programming problem of joint capacity expansion and unit commitment. 
The authors proposed using HAC with Ward's linkage as a primary clustering algorithm, while iteratively picking the optimal number of clusters in the following way.
First, the HAC produces the reduced scenario set of size $N$ (equal to 1 on the first iteration). Next, the optimization problem is solved with the reduced set to produce a lower bound (LB) for the optimal solution of the original complete problem. After that, the first-stage variables which were obtained with the reduced set are fixed and the scheduling stage is recomputed using the full scenario set, providing an upper bound (UB) for the optimal solution. 
Finally, the optimality gap is calculated, and if it is below a set threshold, the algorithm terminates, otherwise, it proceeds onto the next iteration and increases the reduced scenario set size by a step $N^* = N + \text{step}$, where $N^*$ is the new number of clusters. 
The authors found that this algorithm performs better than the method of estimating the optimal number of clusters based on the minimization of the percentage change in within-cluster variance, as the latter tends to overestimate the number of representative periods required to obtain the desired accuracy.

Gonzalez-Castellanos and Pozo \cite{GonzalezCastellanos2021} used the generalized adaptive partition-based method (GAPM) \cite{RamirezPico2021} on the two-stage stochastic problem of optimal wind investment. 
On the first iteration, the optimal partition is equal to the whole scenario set, while the lower and upper bounds are set to negative and positive infinity, respectively. 
Next, the upper bound is calculated similarly to the method discussed above, or if that is not possible, then the current partition is checked against a set of optimality conditions. 
If the gap between the UB and LB is lower than the set constant or the partition satisfies the optimality criteria, then the algorithm stops. 
Otherwise, it proceeds to evaluate the second-stage subproblem for each partition, then moves into the disaggregation stage, where it produces the new partition and proceeds to the next iteration. 
The disaggregation procedure is performed by partitioning the scenarios based on their dual values, i.e., scenarios from the same partition belong to the same subpartition if the value of their dual multipliers is the same.
The authors compared the method against a number of heuristic, clustering and scenario reduction algorithms and found that the GAPM offers by far the best solution, yielding small optimality and solution gaps  in a computationally efficient manner.

\subsection{Spatial aggregation}
\label{sec: Spatial aggregation}

Spatial aggregation refers to methods used to reduce the ``physical'' size of the ESM, i.e., the number of nodes, connections, and generating units. 
The main objective of these methods is to make the network topology less complex by grouping similar or redundant parts of the network. 
Spatial aggregation has been less popular in the literature compared to TSA, and the techniques are mostly limited to two main methods: heuristic algorithms \cite{Palmintier2011, Palmintier2014, Shayesteh2016} and spectral clustering \cite{SanchezGarcia2014, Vankov2020}.

Two studies \cite{Palmintier2011, Palmintier2014} employed simple rule-based heuristic aggregation methods to group generation units based on similar characteristics like fuel type and technology, and the first additionally clustered the units by type plus certain additional characteristics (like age) or aggregated on the plant-level, i.e., all units of the same type inside one facility. 
These methods were employed for spatial aggregation in application to the UC problem, and plant-level aggregation was found to reduce accuracy the least \cite{Palmintier2014}.

A more sophisticated heuristic algorithm was employed by Shayesteh et al. \cite{Shayesteh2016}, where the authors used available transfer capability (ATC) between different buses as the partitioning criterion. The reasoning is that the buses with high ATC have large bulk power transfer between them, which leads to price convergence, hence making those buses effectively the same. This approach was used for network aggregation for the OPF and stochastic UC problems. The authors found that for the OPF problem, the method induced larger errors. However, for the UC problem, the method substantially reduced computational time without significantly decreasing accuracy.

Spectral clustering essentially reduces dimensionality by using the eigenvalues of the similarity matrix constructed for the data \cite{vonLuxburg2007}. Objects and their pairwise similarity measures could be represented through a symmetric similarity matrix $A$ with $A_{ij} \geq 0$. Let $\mathcal{D}$ be a diagonal matrix, where $D_{ii} = \sum\limits_j A_{ij}$, then the Laplacian matrix $L = \mathcal{D} - A$. 
The Laplacian matrix could also be defined differently as per \cite{JianboShi2000}, where the authors used a normalized Laplacian matrix defined as $L = I - \mathcal{D}^{-1/2} A \mathcal{D}^{1/2}$. The latter approach is scale-independent, making it preferable for clustering purposes.

The spectral clustering procedure is the following:
\begin{enumerate}
    \item Calculate the Laplacian matrix.
    \item Calculate $k$ eigenvectors corresponding to the $k$ smallest eigenvalues of $L$.
    \item Form a matrix from these $k$ eigenvectors, where row $i$ defines features of the graph node $i$.
    \item Use a clustering algorithm with these features to produce the partitions.
\end{enumerate}

Sanchez-Garcia et al. \cite{SanchezGarcia2014} applied spectral clustering with HAC on the last step to reduce a power grid. Three similarity measures were tested for the construction of matrix $A$: topology ($A_{ij} = 1$ if nodes $i$ and $j$ are connected), admittance ($A_{ij} = Y_{ij}$, where $Y$ is the admittance matrix) and the average power flow ($A_{ij} = \frac{1}{2} \left( |P_{ij}| + |P_{ji}| \right)$).
The authors found that static measures like the topology or admittance are useful for analysis of the internal structure of the network, while dynamic measures like power flow are better for decision-making that depends on the operating conditions.

Another interesting study was done by Vankov et al. \cite{Vankov2020}, where the authors combined spectral clustering for network aggregation together with HAC to aggregate the time-series data inside the nodes, thus combining spatial and temporal aggregation. They used spectral clustering with k-means on the last step, where the weights in the similarity matrix $A$ were populated using Euclidean, MJC or DTW distance measure for the generation-demand time series for the nodes. After that, the authors used a consensus clustering technique \cite{Lancichinetti2012}, where spectral clustering was applied separately for each day $i$ to build the consensus matrix, while HAC was applied over the consensus matrix to produce the final clusters. The authors showed that this spatio-temporal aggregation approach offered better results than spatial-only methods.

\section{Bibliometric analysis}
\label{sec:Bibliometrics}

We conducted a bibliometric analysis for the papers discussed in the previous section and referenced in Table \ref{table:Methods}. In addition, we extracted from Elsevier's Scopus database additional documents that reference our original set of papers to rank the impact of those papers. We refer to the original 30 papers as first generation, while the documents that cite them are labeled as second generation. We then classified all works into communities of knowledge and applications based on networks metrics.

Table \ref{table: References summary} contains the general summary of the papers used for the analysis. First, it is interesting to note that only 24 documents out of 759 are single-authored, while others contain two or more authors, with an average of almost 4 authors per paper. Secondly, the referenced works were cited around 13 times on average. Finally, the works come from a total of 243 different sources - journals, conference proceedings, books, etc., with the majority being either journal papers (articles and review papers) or conference papers.

\begin{table}[htb]
    \centering
    
    \caption{Summary of the literature collection}
    \label{table: References summary}
    
    \begin{tabular}{lr}
        \toprule
        \textbf{Description} & \textbf{Results} \\
        \midrule
        General information & \\
        \midrule
        Timespan & 2009:2022 \\
        Documents & 759 \\
        Sources & 243 \\
        Average citations per document & 12.73 \\
        References & 41,509 \\
        Authors & 2,019 \\
        Single-authored documents & 24 \\
        Average authors per document & 3.95 \\
        \midrule
        Document types & \\
        \midrule
        Article & 567 \\
        Conference paper & 131 \\
        Review & 47 \\
        Book or book chapter & 13 \\
        Editorial & 1 \\
        \bottomrule
    \end{tabular}
\end{table}

\begin{figure}[h!]
    \centering
    \includegraphics[width=\textwidth, trim={5cm 1cm 7cm 0.5cm}, clip]{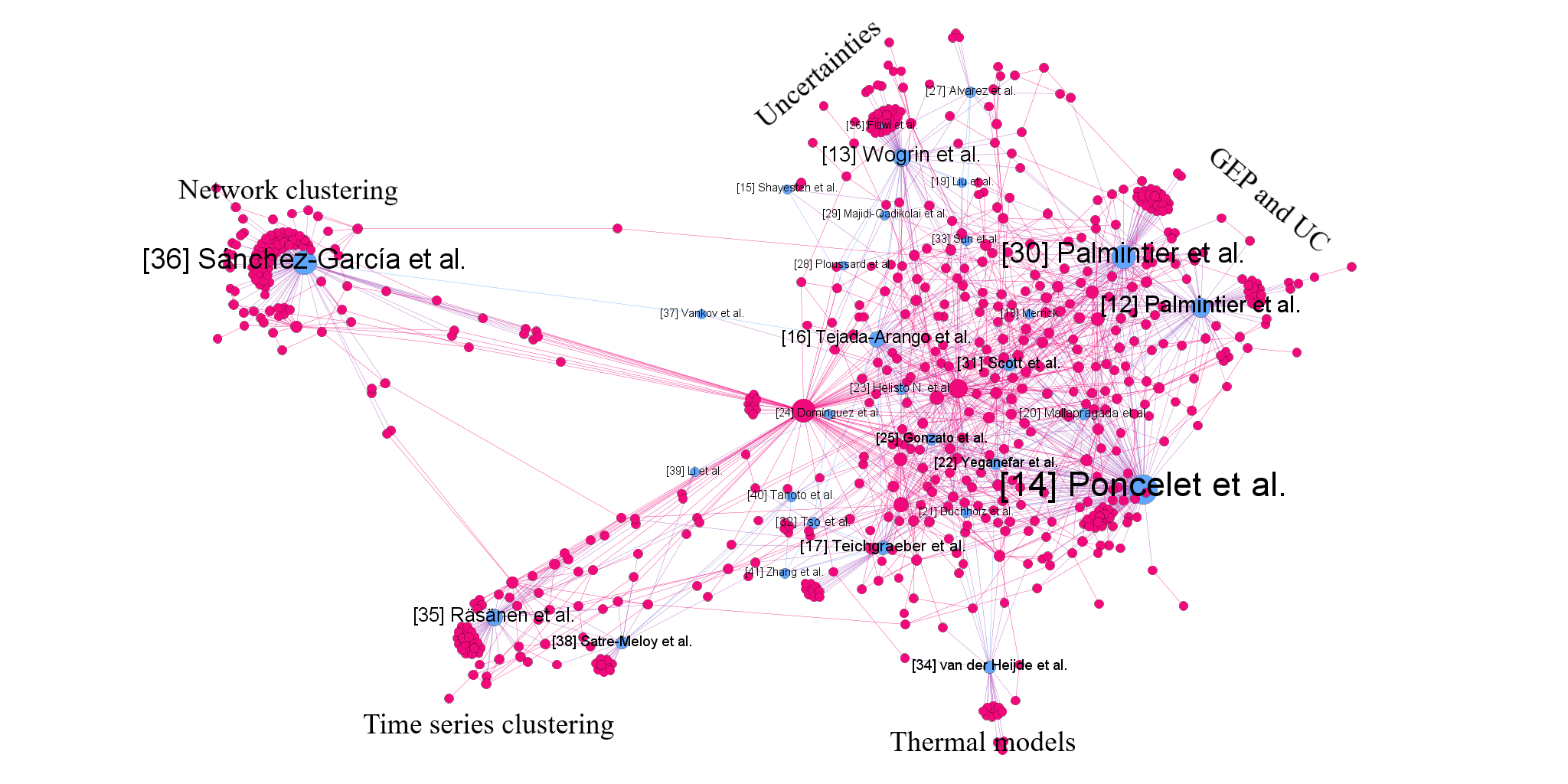}
    \caption{Citation network}
    \label{fig: citation_network}
\end{figure}

Figure \ref{fig: citation_network} shows the citation network built, illustrating the references between works.
Generation-one papers are represented in blue, while generation-two works are in red. 
There are several distinct research groups - papers connected to network clustering (spatial reduction of power grids, graph clustering methods), works focused on TSA (consumer demand clustering, TSA methods, data mining), district heating models, generation expansion planning and unit commitment models, and papers on handling uncertainties (mostly RES). 
The graph was constructed using Gephi \cite{Gephi}, an open source graph visualization software with the Force Atlas layout.

\begin{figure}[h!]
    \centering
    \includegraphics[width=0.8\textwidth]{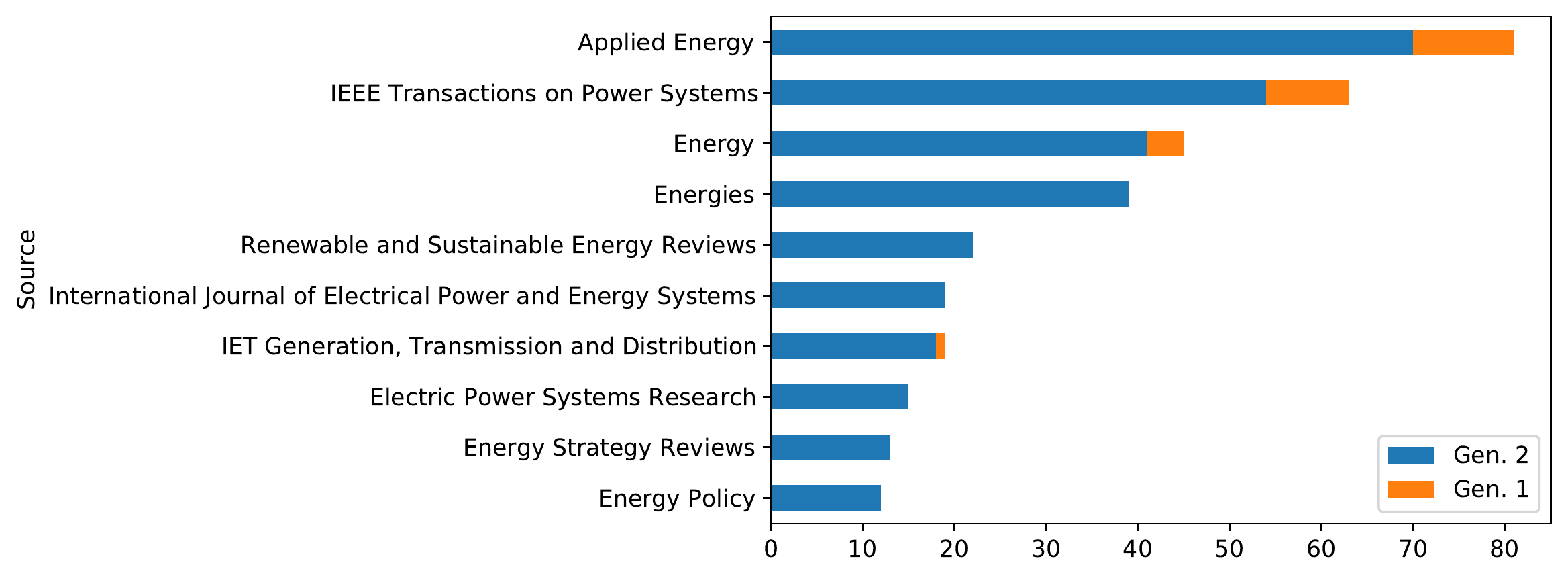}
    \caption{Top 10 journals by publications}
    \label{fig: top10_journals}
\end{figure}

Figure \ref{fig: top10_journals} shows the top 10 journals with the highest number of publications for the papers used in the analysis. Applied Energy and the IEEE Transactions on Power Systems journals are the most popular, with a high number of both Gen. 1 and Gen. 2 papers being published in them.

\begin{figure}[h!]
    \centering
    \includegraphics[width=0.6\textwidth]{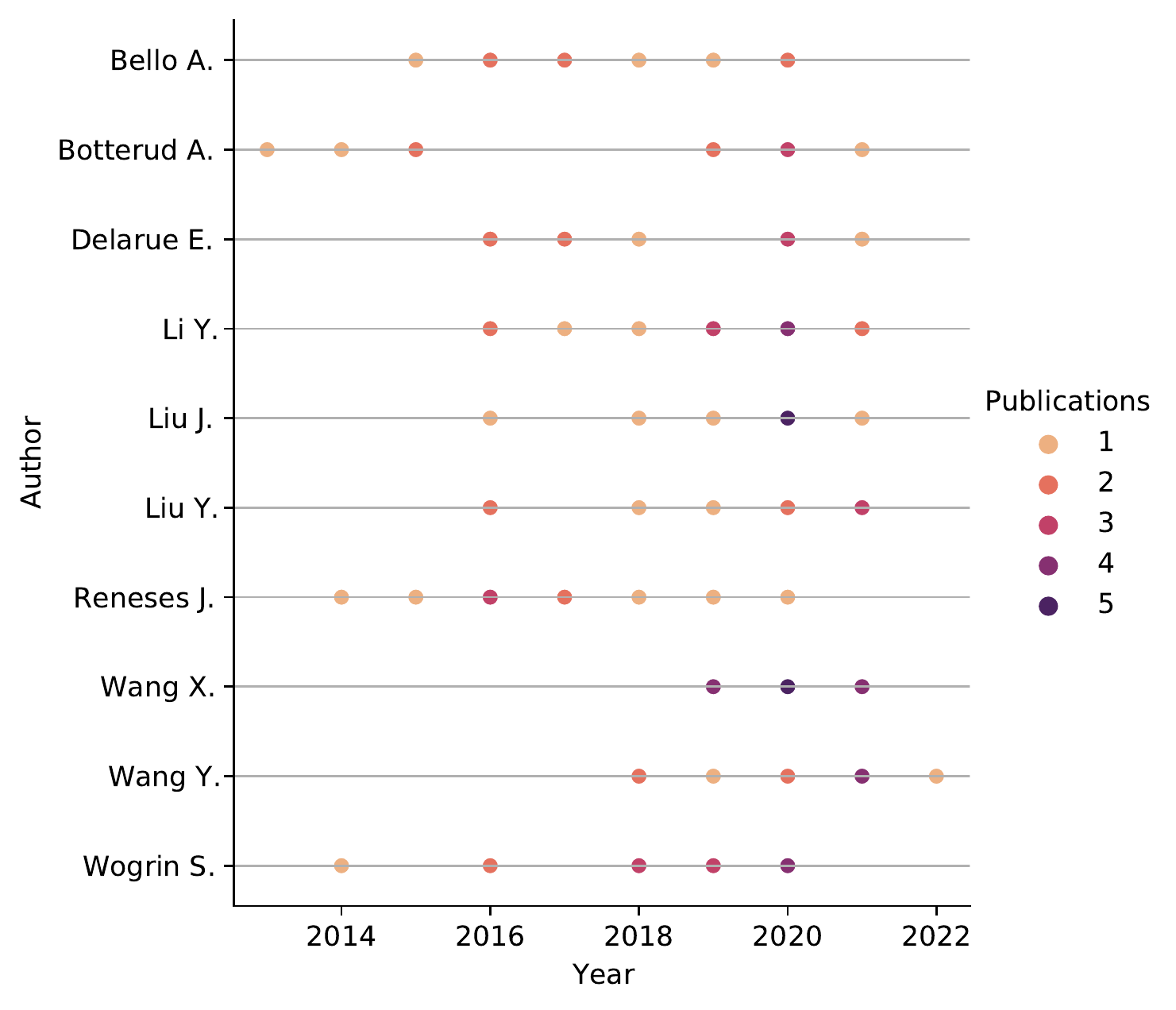}
    \caption{Top 10 authors by publications and their publications per year ordered by family name}
    \label{fig: top10_authors}
\end{figure}

Finally, Figure \ref{fig: top10_authors} contains the top 10 most published authors. For each author, we also plotted the number of published/in press papers per year for the period from 2013 to 2022.

\section{Discussion and recommendations: numerical experiments for the transmission expansion planning problem. }
\label{sec:Experiments}

Determining what constitutes an adequate aggregation methodology for ESM is not an easy task. In the majority of cases, it is not possible to give an a priori response to this question. The choice of an aggregation method is linked to the mathematical model of application and input data. Several works have attempted to benchmark different aggregation methods, but so far there is no consensus about the most competent aggregation method. On the contrary, several works have highlighted the impossibility of providing direct recommendations for selecting the best aggregation method.

This section benchmarks 16 of the most used aggregation methods in the transmission expansion planning problem and two case studies. Performance validation metrics are discussed and used for our analysis. We note that no aggregation method outperforms others in all metrics. It highlights the impossibility of giving, a priori, a solution for the \textit{best} method. However, we note that it is possible to prescribe \textit{adequate} aggregation methods that would satisfy some properties, such as estimating the optimality error bound.

\subsection{Methodology}

Numerical experiments were performed on a transmission expansion planning (TEP) model, as it is one of the recurring and important problems in energy systems literature.
We performed only the temporal aggregation, as the primary objective was to test commonly used algorithms on a single problem to assess their performance, and temporal aggregation has a larger diversity of methods compared to spatial aggregation, where mostly spectral clustering or heuristic methods are used. 
The experiments were performed on two transmission networks, the Garver 6-bus and HRP 38-bus systems, to test the methods' performance on both the smaller and larger systems with different levels of RES penetration.

The following TSA methods were considered: k-means, k-medoids, HAC with single, average, complete and Ward's linkage criteria, FSA with 3 cost functions, and \textit{maximum dissimilarity algorithm} (MDA) \cite{Camus2011}. 
The latter algorithm can be regarded as a scenario reduction algorithm, because it aims to build a reduced scenario set while maximizing the dissimilarity between scenarios in that set, seeking to build the most diversified set. 
The mathematical representation of MDA can be found in \cite{Camus2011}. To our knowledge, it has not been applied for scenario reduction in ESMs before.
Additionally, we used centroid and medoid representations for clustering algorithms, leading to a total number of 16 methods and their variations.

We adopted the following methodology for the numerical experiments:
\begin{enumerate}
    \item The problem is computed using the whole scenario set $\Omega$ to obtain the optimal investment decisions and objective value as reference.
    \item Each reduction algorithm is run $K$ times to produce $K$ reduced sets with different sizes.
    \item The TEP problem is solved for the obtained reduced sets, then we choose the reduced set $\Omega_R$ which leads to the lowest objective value gap when compared with the full set solution and treat the corresponding reduced set size $|\Omega_R|$ as optimal for this algorithm.
    \item We run an out-of-sample (OOS) test for the algorithms, where the first-stage decisions are made with the reduced set, while the second-stage is evaluated using these fixed decisions with the disaggregated full scenario set. This way, we can assess how well the reduced set $\Omega_R$ approximates the original set $\Omega$.
    \item We compare the methods using the optimal reduced set size obtained in step 3 and its corresponding out-of-sample objective function gap from step 4.
    \item We compare the run time of each algorithm with the highest reduced set size.
\end{enumerate}

The number of reduced set sizes $K \in \{1, \dots, 10\}$. 
For the Garver network, we also performed a sensitivity analysis to gauge how variation in model parameters affects computational time, highlighting the importance of using TSA to reduce the computational burden.

The numerical experiments were performed on a cluster with 256 GB of RAM and 2x Intel Xeon Gold 6148 processors. The code was written in the Python 3 programming language, using Gurobi as a solver and with the MIP gap set to 0.1\% \cite{Gurobi}.

\subsection{Transmission expansion planning modeling}

The model formulation for the employed transmission expansion planning problem is presented below \eqref{eq: tep_obj}--\eqref{eq: tep_binary}, and its notation is given in Table \ref{table: Notation}:

\begin{table}[h!]
    \centering
    
    \caption{Notation for the transmission expansion planning (TEP) problem}
    \label{table: Notation}
    
    \resizebox{\textwidth}{!}{
    \begin{tabular}[t]{ccc}
        \toprule
        \textbf{Sets} & \textbf{Variables} & \textbf{Parameters} \\
        \midrule
        \makecell[lc]{$ \Omega_T $ - time periods \\ $ \Omega_K $ - scenarios \\ $ \Omega_{Le} $ - existing lines \\ $ \Omega_{Lc} $ - candidate lines \\ $ \Omega_L $ - all lines \\ $ \Omega_G $ - conventional units \\ $ \Omega_{Rg} $ - renewable units \\ $ \Omega_G^n $ - conventional units at node $n$ \\ $ \Omega_{Rg}^n $ - renewable units at node $n$ \\ $ \Omega_D^n $ - demand at node $n$}
        & 
        \makecell[lc]{$ x_l \in \{0,1\} $ - indicator showing whether \\ line $l$ exists \\ $ P_g^{t} \in \mathbb{R}  $ - power produced by the gene-\\rating unit $g$ at time $t$ \\ $ f_l^t \in \mathbb{R}  $ - power flow through line $l$ \\ at time $t$ \\ $ \theta_n^t \in \mathbb{R} $ - voltage phase angle at node $n$ \\ at time $t$ \\ $ f_l^{max} $ - current capacity of line $l$}
        &
        \makecell[lc]{$ c_l^{inv} \in \mathbb{R} $ - investment cost of line $l$ \\ $ c_g^{op} \in \mathbb{R} $ - operational cost of generator $g$ \\ $ B_l \in \mathbb{R} $ - susceptance of line $l$ \\ $ F_l^{max} $ - maximum capacity of line $l$ \\ $ P_g^{max} $ - maximum capacity of \\ conventional unit $g$ \\ $ R_r^t $ - maximum capacity of \\ renewable unit $r$ at time $t$ \\ $ \Delta P_{\text{DW}}^g $ - down ramping limits for unit $g$ \\ $ \Delta P_{\text{UP}}^g $ - up ramping limits for unit $g$}\\
        \bottomrule
    \end{tabular}}
\end{table}

\begin{alignat}{2}
    \min & \sum\limits_{l \in \Omega_{Lc}} \left[ c_l^{\text{inv,fix}} x_l + c_l^{\text{inv,var}} f_l^{max} \right] + \sum\limits_{k \in \Omega_K} \left[ \sum\limits_{g \in \Omega_g}
    \sum\limits_{t \in \Omega_T} c_{gk}^{op} P_{gk}^t \right] \label{eq: tep_obj}
    \\
    \text{s.t. } & \sum\limits_{g \in \Omega_g^n} P_{gk}^t + \sum\limits_{r \in \Omega_{Rg}^n} P_{rk}^t - \sum\limits_{l | i(l) = n} f_{lk}^t + \sum\limits_{l | j(l) = n} f_{lk}^t = \sum\limits_{d \in \Omega_D^n} D_{dk}^t \quad && t \in \Omega_T, l \in \Omega_L, k \in \Omega_K \label{eq: tep_pb}
    \\
    & x_l F_{l}^{min} \leq f_l^{max} \leq x_l F_{l}^{max} \quad && l \in \Omega_L \label{eq: tep_fmax}
    \\
    & -f_{l}^{max} \leq f_{lk}^t \leq f_{l}^{max} \quad && t \in \Omega_T, l \in \Omega_L, k \in \Omega_K \label{eq: tep_ft_max}
    \\
    & -(1 - x_l) M \leq f_{lk}^t + B_l \left( \theta_{i(l) k}^t - \theta_{j(l) k}^t \right) \leq (1 - x_l) M \quad && l \in \Omega_L, k \in \Omega_K, t \in \Omega_T \label{eq: tep_pf}
    \\
    & -\pi \leq \theta_{i(l) k}^t - \theta_{j(l) k}^t \leq \pi \quad && l \in \Omega_L, k \in \Omega_K, t \in \Omega_T \label{eq: tep_angle}
    \\
    & 0 \leq P_{gk}^t \leq P_{g}^{max} \quad && g \in \Omega_G, k \in \Omega_K, t \in \Omega_T \label{eq: tep_cg}
    \\
    & 0 \leq P_{rk}^t \leq R_{rk}^t \quad && k \in \Omega_K, r \in \Omega_{Rg}, t \in \Omega_T \label{eq: tep_rg}
    \\
     & \Delta P_{\text{DW}}^g \leq P_{gk}^t - P_{gk}^{t-1} \leq \Delta P_{\text{UP}}^g \quad && k \in \Omega_K, g \in \Omega_G, t \in \Omega_T \label{eq: tep_ramp}
    \\
    & x_l = 1 \quad && l \in \Omega_{Le} \label{eq: tep_le}
    \\
    & x_l \in \{0, 1\} \quad && l \in \Omega_{Lc} \label{eq: tep_binary}
\end{alignat}

The objective function \eqref{eq: tep_obj} consists of two parts: the first summation refers to the cost of investment in building new lines, while the second refers to operational expenses. 
The operational expenses are essentially the same as in other two-stage stochastic TEP models. The novel part of the proposed MILP model is the line investment part of the objective function, which is driven by two costs: one is $c_l^{\text{inv,fix}}$, which represents the cost of building the essential line infrastructure as well as some minimal working capacity, multiplied by the indicator variable $x_l$, which shows whether line $l$ was built or not; the second is $c_l^{\text{inv,var}}$, which is the cost of line capacity expansion, multiplied by the current maximum line capacity $f_l^{max}$. Such a formulation is closer to reality, as it both represents the ``lumpy'' characteristics of the transmission investments and still has a reasonable level of control over the line capacity, rather than representing the investment decision solely by binary indicators (like in most MILP models) or discarding the associated building costs (like in LP models).

Constraint \eqref{eq: tep_pb} represents the nodal power balance. Constraints \eqref{eq: tep_fmax}-\eqref{eq: tep_angle} are the power flow constraints. The minimum and maximum line capacity is set by \eqref{eq: tep_fmax}.
Bounds \eqref{eq: tep_ft_max} indicate that the power flow through the line cannot exceed its capacity.
The line flow is given by \eqref{eq: tep_pf}, while \eqref{eq: tep_angle} represents the nodal angle limits. Another set of constraints \eqref{eq: tep_cg}-\eqref{eq: tep_ramp} is related to the system's the generators. Expression \eqref{eq: tep_cg} represents the conventional generation limits, \eqref{eq: tep_rg} represents the renewable generation limits and depends on the scenario realization, and \eqref{eq: tep_ramp} is the ramping constraint for conventional generators. 
Finally, the last two constraints, \eqref{eq: tep_le} and \eqref{eq: tep_binary}, determine whether the line is built or not and the binary nature of the decision variable $x_l$.

\subsection{The Garver 6-bus system}

\subsubsection{Data}

The Garver network was first introduced in \cite{Garver1970} and has been frequently used as an illustrative case for TEP studies in the literature. 
We modified the network by replacing the conventional thermal generator at node 1 with a wind generator of equivalent capacity. 
The candidate line set consists of unbuilt lines 1-3, 1-6, 2-5, 2-6, 3-4, 3-6, 4-5, 4-6, 5-6, and all the existing lines, which are considered for capacity expansion. 
The minimum capacity required to justify building one of the unbuilt lines is set to 20\% of its maximum value.
The maximum number of lines in one corridor is three, which is equivalent to expanding the lines up to 300\% of their maximum capacity. 
Hourly wind and load scenarios for one year were taken from the ENTSO-E transparency data set for Germany \cite{ENTSOE}. 


\subsubsection{Loadability sensitivity analysis}

\begin{figure}[h]
    \centering
    \begin{subfigure}{0.48\textwidth}
     \includegraphics[width=\linewidth]{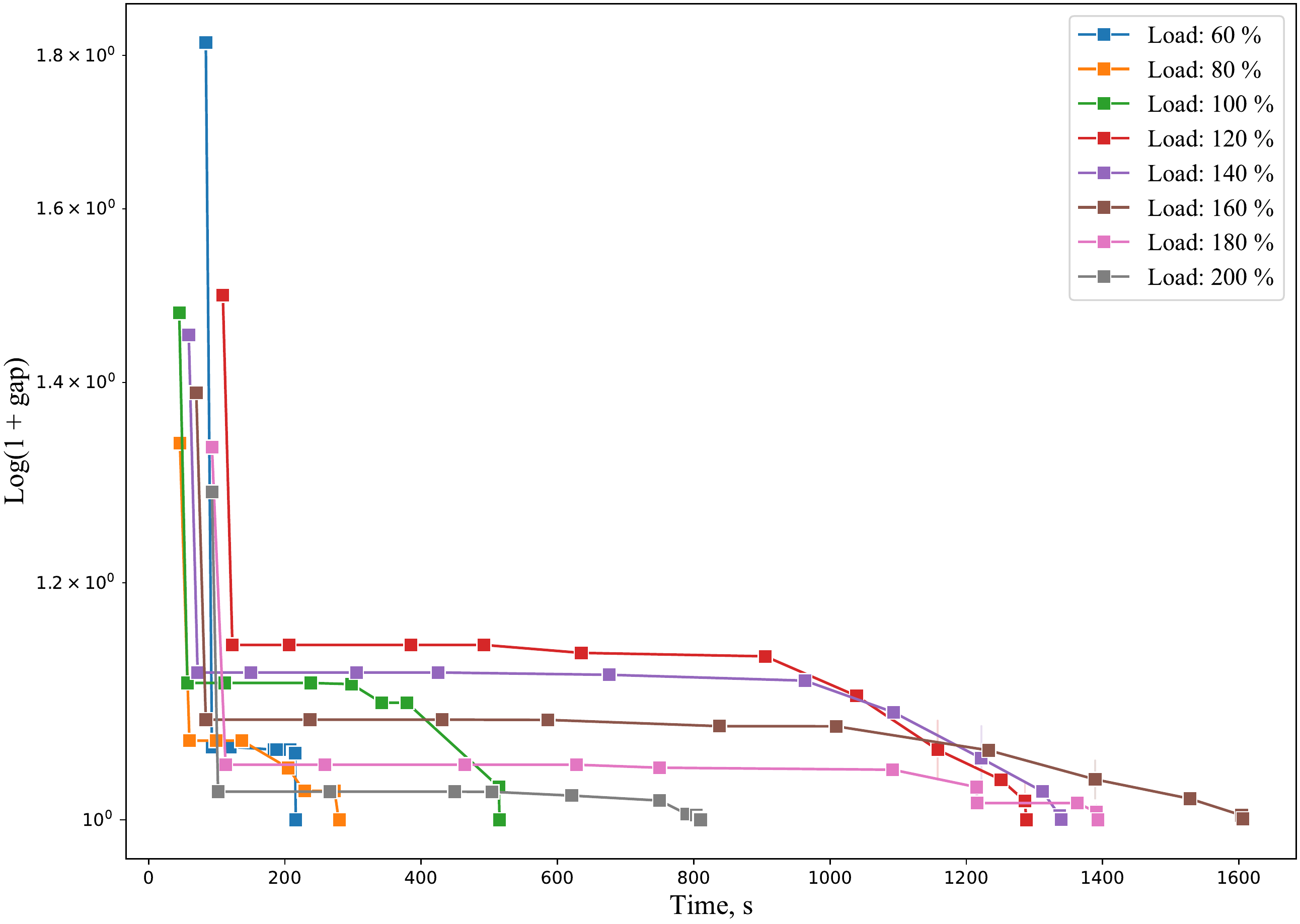}
        \caption{MIP gap evolution in time with different load levels}
                \label{fig: BB_Exploration_MIP}
    \end{subfigure}%
    \hspace{0.05cm}
    \begin{subfigure}{0.48\textwidth}
        \includegraphics[width=0.97\linewidth]{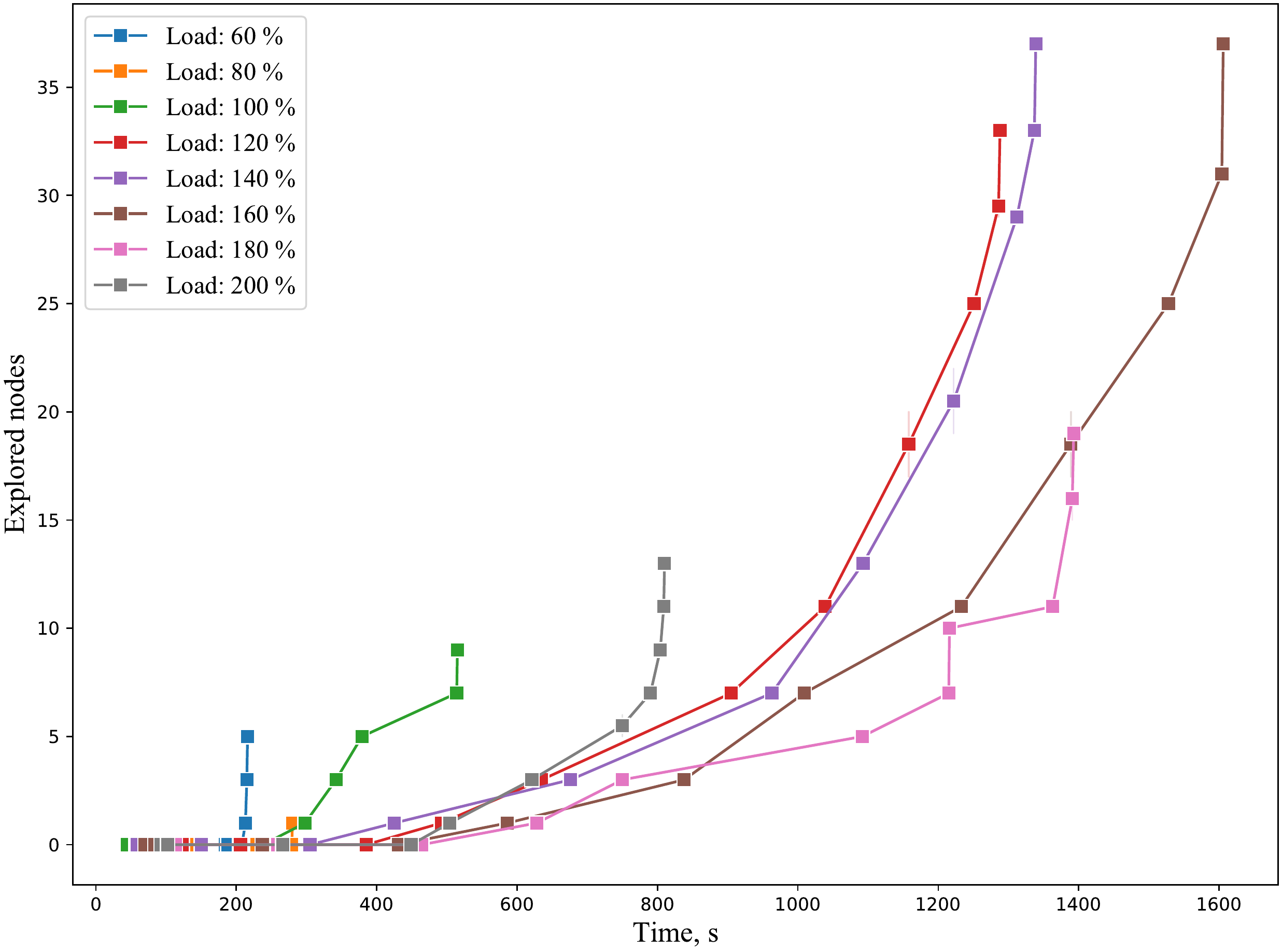}
    \caption{Node exploration progression with different load levels}
                    \label{fig: BB_Exploration_Node}
    \end{subfigure}
    \caption{Branch-and-bound solver exploration progress with different load levels for the Garver 6-bus network}
    \label{fig: BB_Exploration}
\end{figure}

Before comparing the TSA methods, we first ran a sensitivity analysis on this illustrative case with the full scenario set to investigate how varying one of the problem parameters would affect the computational complexity of the problem. 
For this purpose, we varied the demands at each node from 60\% to 200\% of the initial values with a step of 20\%. 
Figure \ref{fig: BB_Exploration_MIP} shows the evolution of the MIP gap (logarithmic scale) with time for the different load scenarios. 
It can be seen that the computational time increases as nodal load levels increase until it reaches the 160\% mark, at which point it starts dropping back.

This behavior can be explained by the NP-hard nature of the optimization problem and the working principle of the algorithm employed to solve such problems, namely the branch-and-bound (B\&B) algorithm \cite{Land1960}.
Briefly, the B\&B algorithm first solves an LP relaxation of the MILP problem and checks whether the LP solution satisfies the integrality constraints. 
If it does so, then the solution has been found; otherwise, the algorithm builds a search tree by splitting the integer variables and builds simpler MIP subproblems that need to be solved. 
Solving these subproblems one by one and picking the optimal integer values for the variables leads to an optimal solution for the original set. 
Therefore, the larger the search tree grows, the more MIP subproblems need to be computed, which in turn drastically increases computational time. 
This means that low load levels can be easily solved as most probably there is no need to build or expand many lines to satisfy consumption, while really high load levels instead prompt almost all lines to be built and expanded. 
Both cases are ``straightforward" in the sense that they do not require the search tree be grown too large or every branching case be considered, while the cases in the middle might require the tree to be exhaustively traversed to find the optimal solution by evaluating many feasible solutions.

This assessment is supported by Figure \ref{fig: BB_Exploration_Node}, which plots the number of explored nodes in the branching tree against time. As can be seen, the number of explored nodes grows with the increasing load percentage, but after peaking at 140\% and 160\%, it drops for the 180\% and 200\% cases. 


\subsubsection{Comparison of methods}

\begin{figure}[h]
    \centering
    \begin{subfigure}{0.32\textwidth}
        \includegraphics[width=\linewidth]{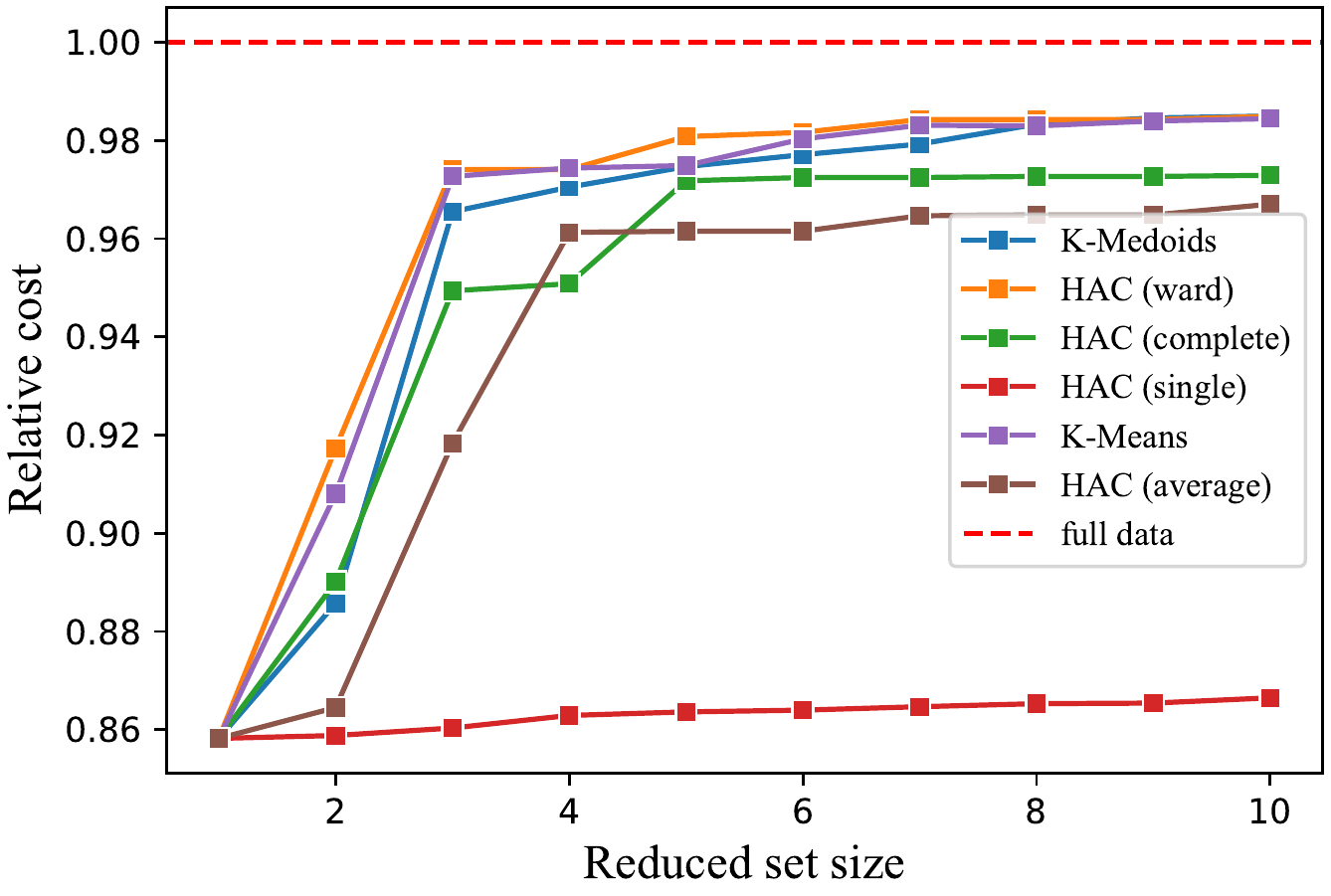}
        \caption{Centroid representation}
    \end{subfigure}%
        \hspace{0.05cm}
    \begin{subfigure}{0.32\textwidth}
        \includegraphics[width=\linewidth]{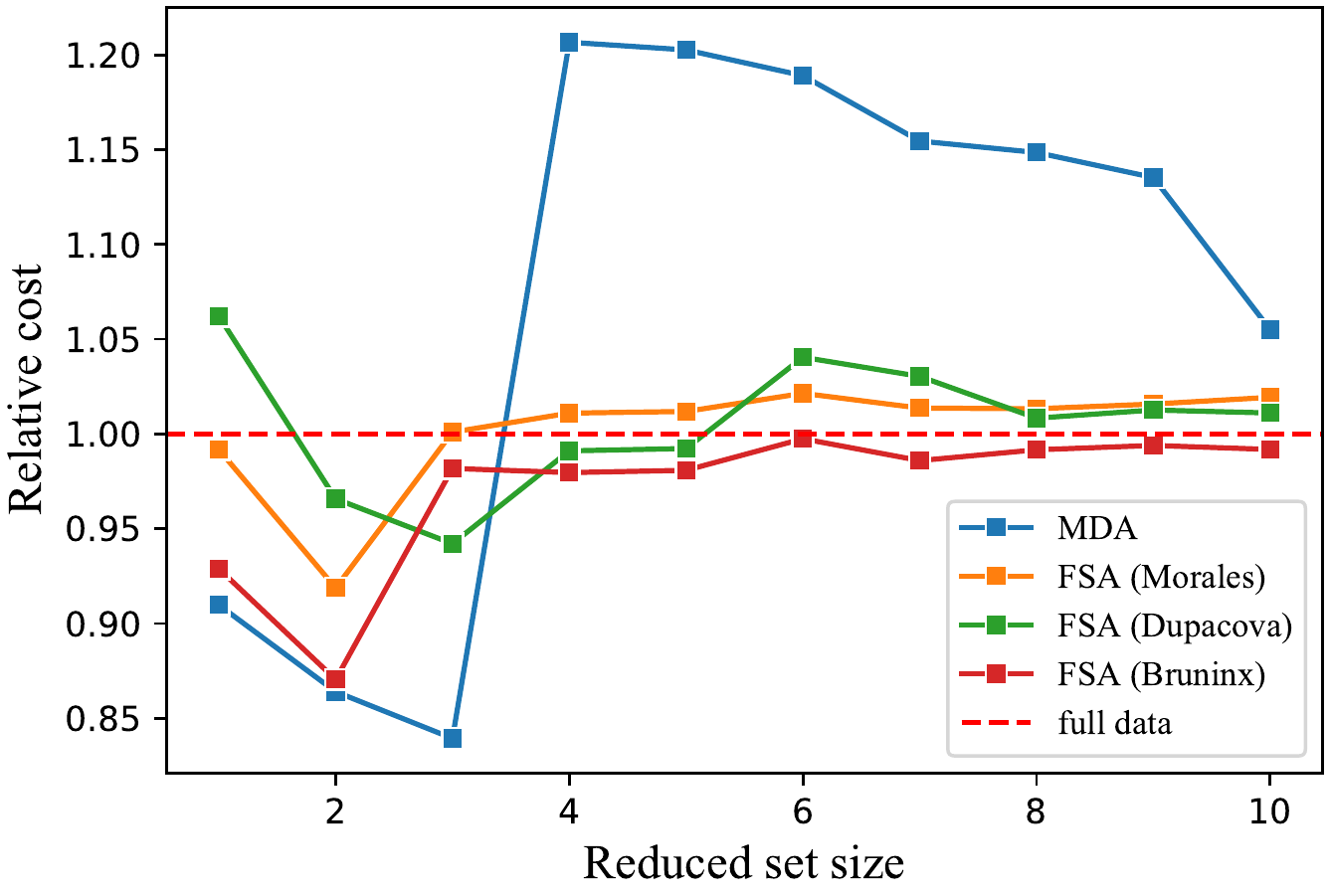}
        \caption{Scenario representation}
    \end{subfigure}
        \hspace{0.05cm}
    \begin{subfigure}{0.32\textwidth}
        \includegraphics[width=\linewidth]{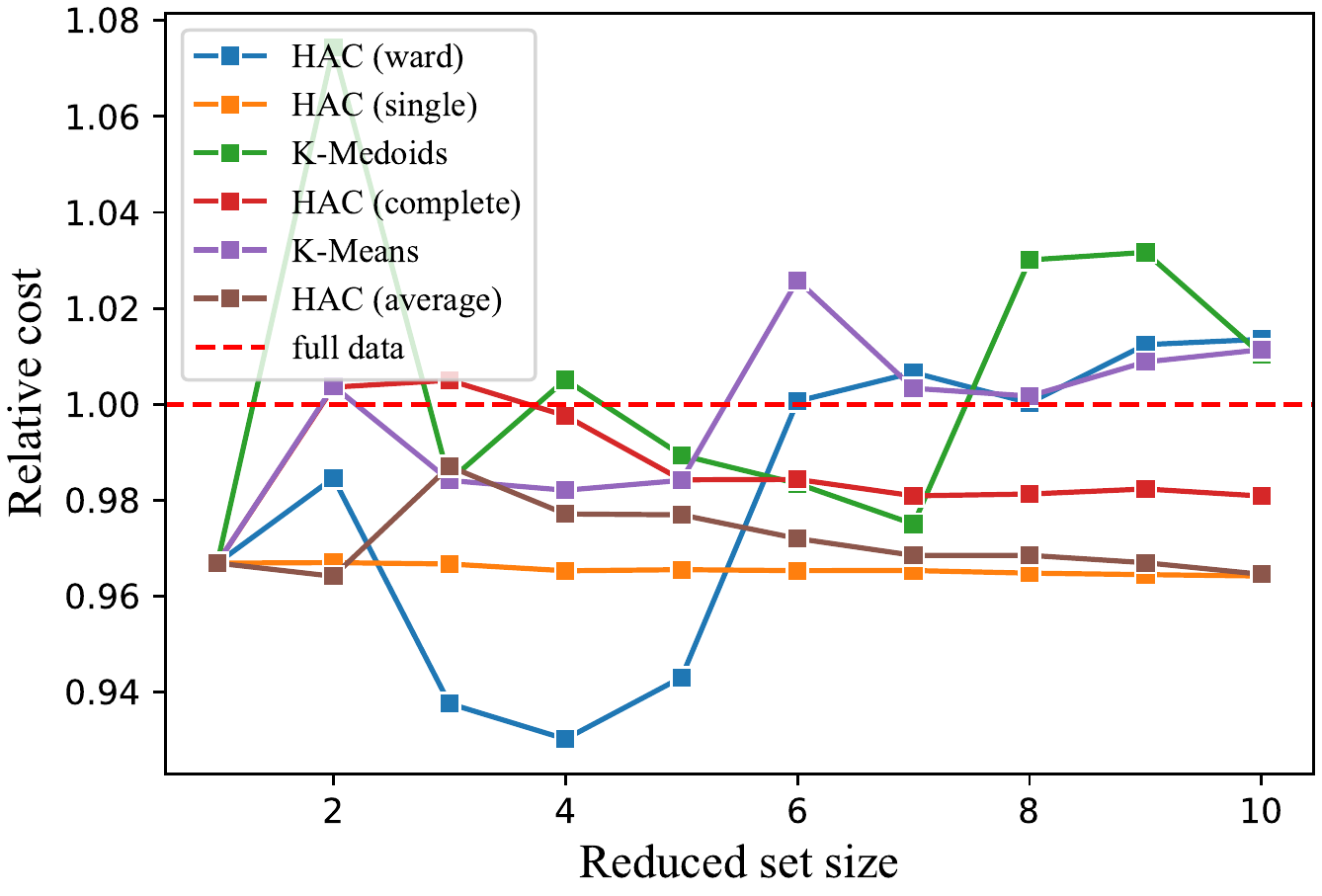}
        \caption{Medoid representation}
    \end{subfigure}
    \caption{In-sample relative objective value change with different reduced set sizes for the Garver 6-bus network}
    \label{fig: garver_insample}
\end{figure}

Figure \ref{fig: garver_insample} presents the evolution of the in-sample relative objective cost value and the objective value obtained with the reduced set divided by the optimal objective value with all scenarios, with change in the reduced set size for each method. 
In the figure, the scenario representation refers to the reduced set obtained by methods that directly build the scenario set without using representation functions.

Based on the obtained results, we can make the same observation as Teichgraeber et al. \cite{Teichgraeber2019}: the methods with centroid representation produce reduced scenarios sets which can lead to lower bounds for the optimal solution. Additionally, the solutions which use HAC with single linkage do not change very much as the reduced set size increases regardless of the representation function used. 
This is explained by the nature of the single linkage criterion, which tends to amplify \textit{the-rich-get-richer} behavior, resulting in worse diversity of the resulting reduced scenario set and leading to worse performance in the optimization problem. 
Finally, the MDA method leads to a large jump between the reduced set sizes of 3 and 4.
This anomaly stems from the approach used to initialize the scenario set, where we used a scenario with the peak daily load as per \cite{Camus2011}. The peak load day can coincide with the peak generation day, which would result in an underestimation of the objective value, but as the algorithm adds additional scenarios while maximizing their dissimilarity, it will sooner or later move to the expected overestimation of the objective value.

\begin{figure}[h]
    \centering
    \begin{subfigure}{0.32\textwidth}
        \includegraphics[width=\linewidth]{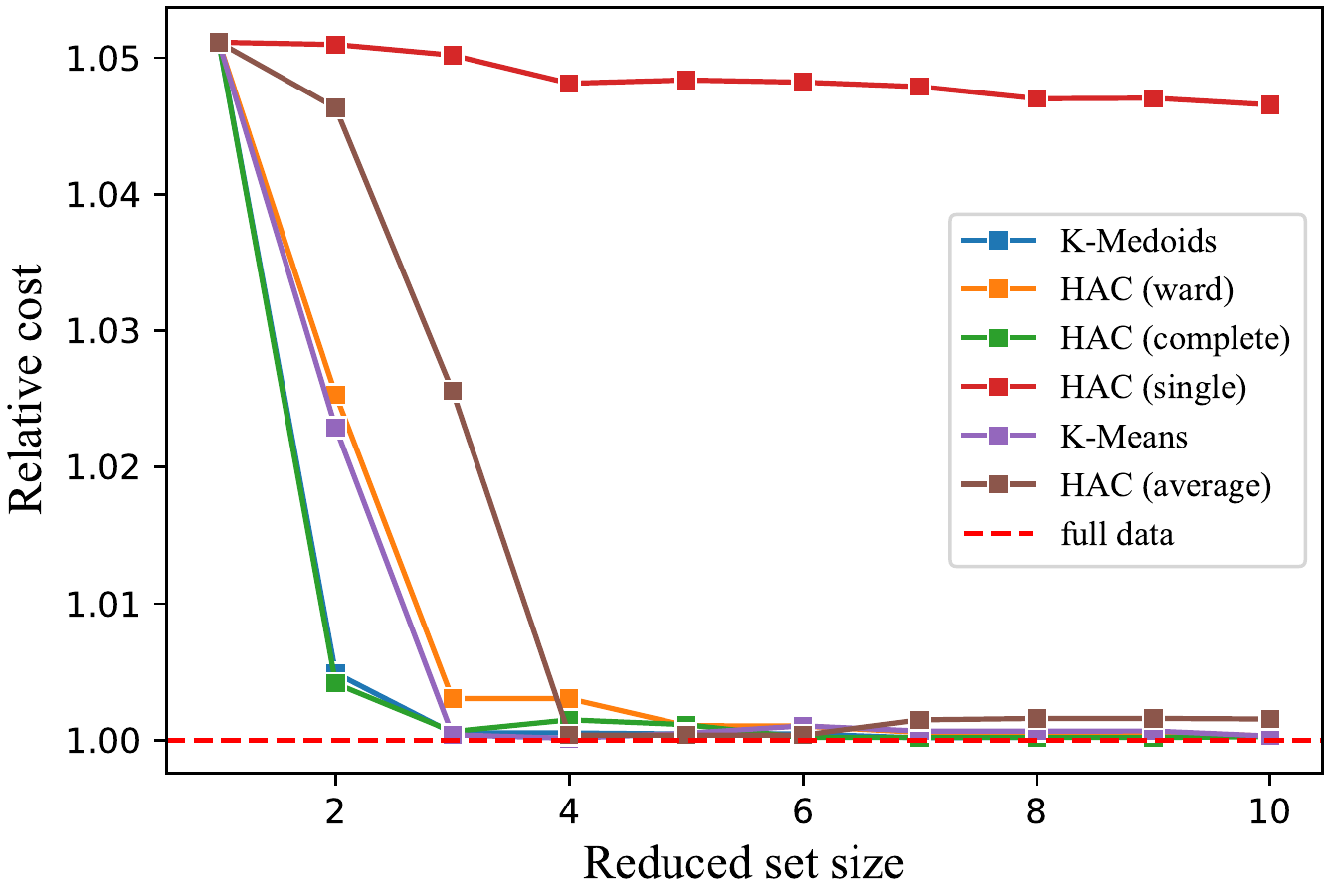}
        \caption{Centroid representation}
    \end{subfigure}%
    \hspace{0.05cm}
    \begin{subfigure}{0.32\textwidth}
        \includegraphics[width=\linewidth]{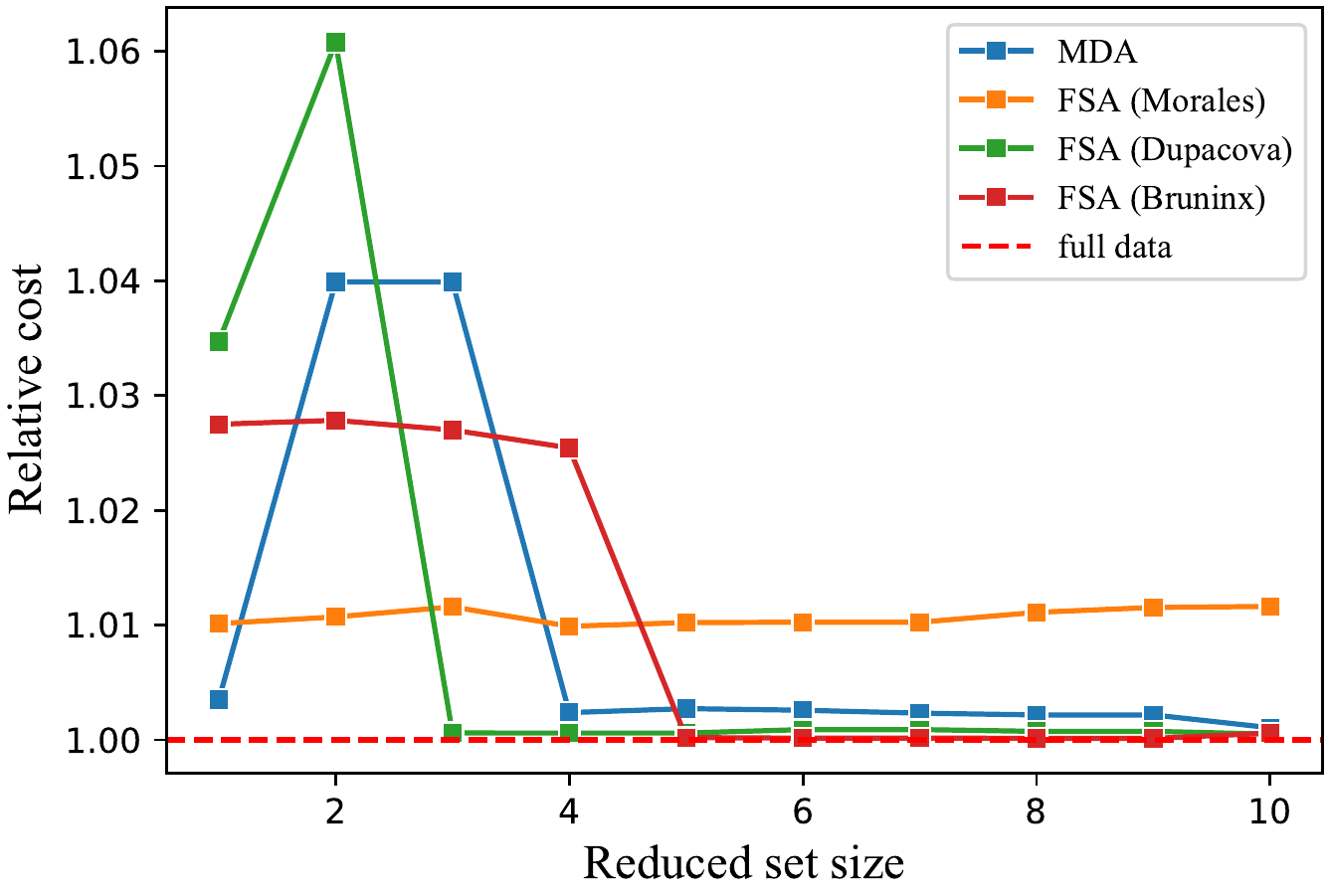}
        \caption{Scenario representation}
    \end{subfigure}
    	\hspace{0.05cm}
    \begin{subfigure}{0.32\textwidth}
        \includegraphics[width=\linewidth]{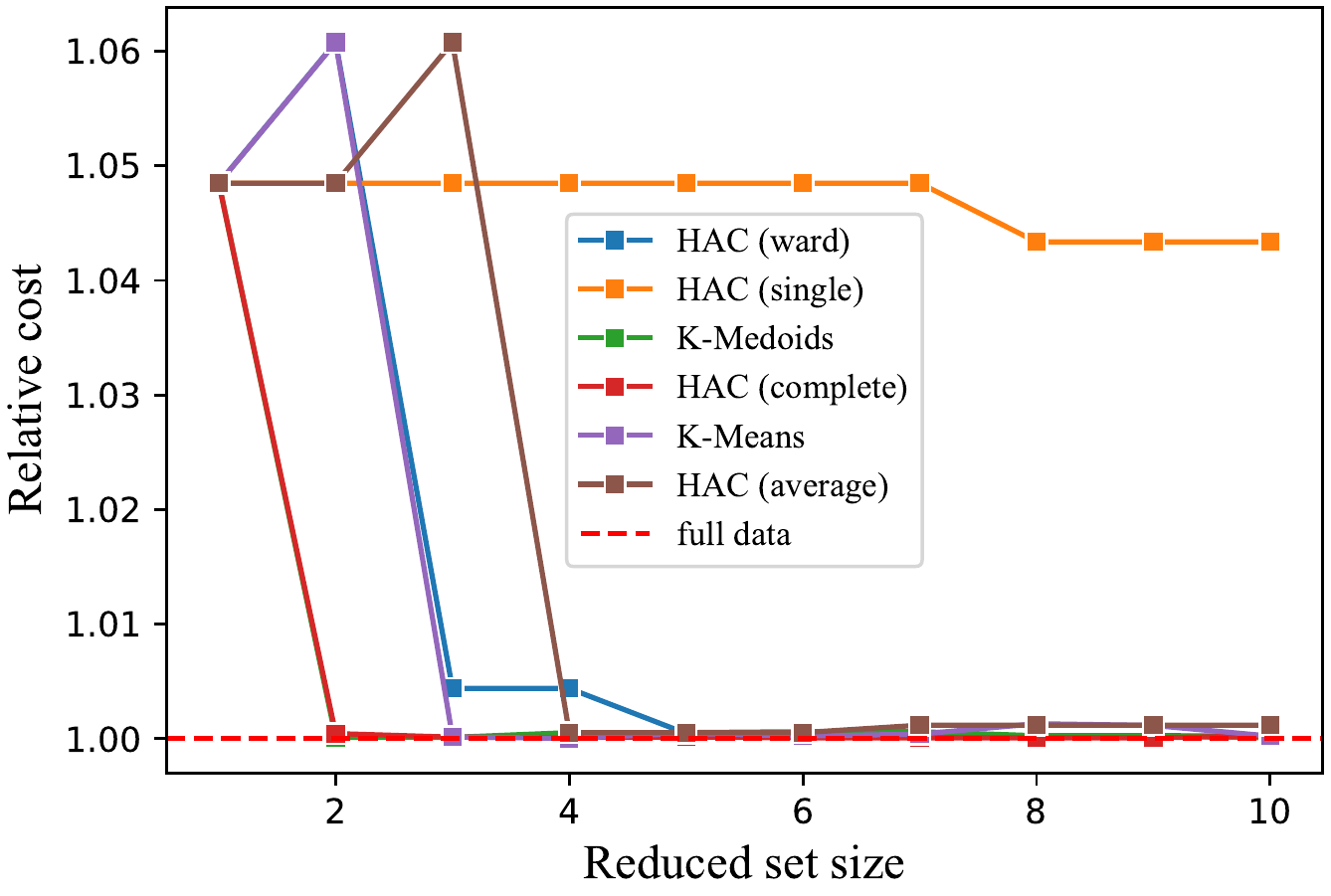}
        \caption{Medoid representation}
    \end{subfigure}
    \caption{Out-of-sample relative objective value change with different reduced set sizes for the Garver 6-bus network}
    \label{fig: garver_oos}
\end{figure}

The out-of-sample relative objective costs can be seen on Figure \ref{fig: garver_oos}.
The results are consistent with the above plots, but it is important to note the performance of the MDA method: the overestimation that was present in the in-sample costs does not translate into the out-of-sample costs, and the diversified scenario set produced by the algorithm approximates the objective value really well.
The first thing to notice is that the behavior of the HAC with single linkage continues on the out-of-sample test, which further confirms that it is connected to the method's process of selection. 
Unexpectedly, the FSA with the cost function from \cite{Morales2009} also produces scenario sets that do not fully capture the diversity of the original set, as is evident from the almost non-changing relative cost in the out-of-sample test, which suggests that using objective values obtained from solving the problem with a deterministic fixed first stage in the cost function might lead to less diversified scenario sets.

Next, we followed the methodology and determined the optimal reduced set size for each algorithm, then compared the respective produced OOS objective gaps and investment decisions. Results are presented in Table \ref{table: garver_methods_comparison}, for scale - 1\% of the total objective gap is about 317,000 USD. The investment part constitutes about 10\% of the total cost, with the rest being associated with operational costs. Only the solution using HAC with average linkage and medoid representation was unable to make correct investment decisions, which led to a large OOS gap. HAC with single linkage, regardless of representation, also led to a rather large OOS gap - while the investments themselves were made in the correct lines, it chose a smaller capacity for line 1-5, which was insufficient and resulted in extensive use of load shedding and increased the end cost. In table, the best performing algorithms in each group are marked in bold. The overall best performing algorithm was HAC with complete linkage, which produced solutions with the best results for both centroid and medoid cluster representations. This means that the algorithm obtained very good reduced sets that were able to capture the majority of the information contained in the original scenario set. 
\begin{table}[h]
    \centering
        \caption{Comparison of algorithms on the Garver 6-bus system}
    \begin{tabular}{ccccc}
        \toprule
        \textbf{Method} & \textbf{OOS gap, \%} & \textbf{Built lines} & \textbf{Capacities, MW} & \textbf{Num. repr.} \\
        \midrule
        All scenarios & - & \{1-5, 2-6, 4-6\} & \{92.4, 300, 247.6\} & 365 \\
        \midrule
        \multicolumn{5}{c}{Centroid representation} \\
        \midrule
        K-Means & 0.03 & \{1-5, 2-6, 4-6\} & \{89.8, 300, 247.1\} & 10 \\
        K-Medoids & 0.05 & \{1-5, 2-6, 4-6\} & \{89.8, 300, 247\} & 10 \\
        HAC (Ward) & 0.05 & \{1-5, 2-6, 4-6\} & \{89.7, 300, 247.5\} & 10 \\
        HAC (average) & 0.15 & \{1-5, 2-6, 4-6\} & \{86.9, 300, 246.1\} & 10 \\
        HAC (single) & 4.65 & \{1-5, 2-6, 4-6\} & \{69, 297, 237.6\} & 10 \\
        HAC (complete) & \textbf{0.02} & \{1-5, 2-6, 4-6\} & \{90.4, 300, 246.4\} & 10 \\
        \midrule
        \multicolumn{5}{c}{Medoid representation} \\
        \midrule
        K-Means & 0.13 & \{1-5, 2-6, 4-6\} & \{89.8, 300, 247.1\} & 8 \\
        K-Medoids & 0.05 & \{1-5, 2-6, 4-6\} & \{89.8, 300, 247\} & 4 \\
        HAC (Ward) & 0.03 & \{1-5, 2-6, 4-6\} & \{90.8, 300, 249.1\} & 8 \\
        HAC (average) & 6.08 & \{2-6, 4-6\} & \{300, 242.3\} & 3 \\
        HAC (single) & 4.85 & \{1-5, 2-6, 4-6)\} & \{65.7, 297, 244.4\} & 2 \\
        HAC (complete) & \textbf{0.01} & \{1-5, 2-6, 4-6\} & \{90.8, 300, 247.5\} & 4 \\
        \midrule
        \multicolumn{5}{c}{Scenario representation} \\
        \midrule
        MDA & 0.10 & \{1-5, 2-6, 4-6\} & \{94.3, 300, 251.1\} & 10 \\
        FSA (Dupačová) & \textbf{0.06} & \{1-5, 2-6, 4-6\} & \{88.8, 300, 248.4\} & 5 \\
        FSA (Morales) & 1.16 & \{1-5, 2-6, 4-6\} & \{78, 300, 248\} & 3 \\
        FSA (Bruninx) & 0.14 & \{1-5, 2-6, 4-6\} & \{90.6, 300, 247\} & 6 \\
        \bottomrule
    \end{tabular}
    
    \label{table: garver_methods_comparison}
\end{table}

Last but not least, we compared the speed of algorithms with $|\Omega_R| = 10$. Results are shown in Table \ref{table: garver_runtime} - solving TEP problem on the Garver system with 10 scenarios takes around 1 second and HAC methods had very similar runtimes regardless of the linkage criterion, so they were aggregated into a single entity. FSA algorithms generally showed slower speeds and longer runtimes, especially the variations which require the computation of different objective values for the cost function. For instance, FSA with cost function \eqref{eq: Bruninx} was extremely slow, while the results obtained in the previous stages were worse than the results from the faster algorithms. This means that using objective values for the cost function in FSA is very expensive computationally but does not necessarily guarantee better results.

\begin{table}[h]
    \centering
        \caption{Running time comparison of algorithms on the Garver 6-bus system}
    \begin{tabular}{ccc}
        \toprule
        \textbf{Method} & \textbf{Method time, s} & \textbf{Total time, s} \\
        \midrule
        All scenarios & - & 509.2 \\
        \midrule
        K-Means & $0.401 \pm 0.007$ & 1.4 \\
        K-Medoids & $0.011 \pm 0.003$ & 1 \\
        HAC & $0.020 \pm 0.001$ & 1 \\
        MDA & $0.020 \pm 0.001$ & 1 \\
        FSA (Dupačová) & $0.743 \pm 0.004$ & 1.6 \\
        FSA (Morales) & $18.2 \pm 0.1$ & 19.1 \\
        FSA (Bruninx) & $305 \pm 1.2$ & 306 \\
        \bottomrule
    \end{tabular}
    \label{table: garver_runtime}
\end{table}

\subsection{High Renewable Penetration (HRP) 38-bus system}


The authors who proposed the High Renewable Penetration (HRP) 38-bus system in \cite{Zhuo2020} were motivated by the lack of standardized test systems for stochastic TEP problems with high levels of renewable penetration. We refer to the original work for the bus, line, and scenario data.



\begin{figure}[h!]
    \centering
    \begin{subfigure}{0.32\textwidth}
        \includegraphics[width=\linewidth]{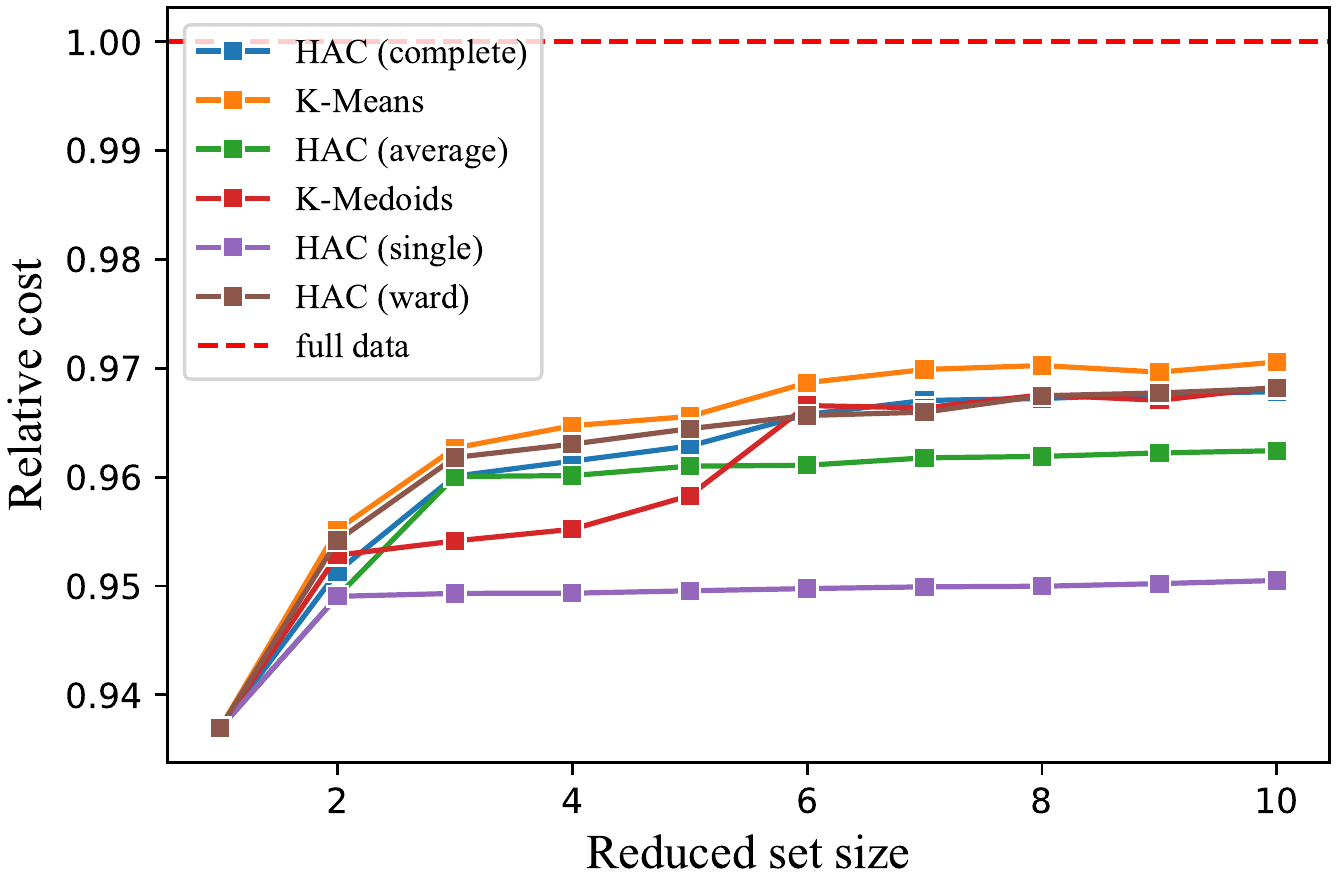}
        \caption{Centroid representation}
    \end{subfigure}%
        \hspace{0.05cm}
    \begin{subfigure}{0.32\textwidth}
        \includegraphics[width=\linewidth]{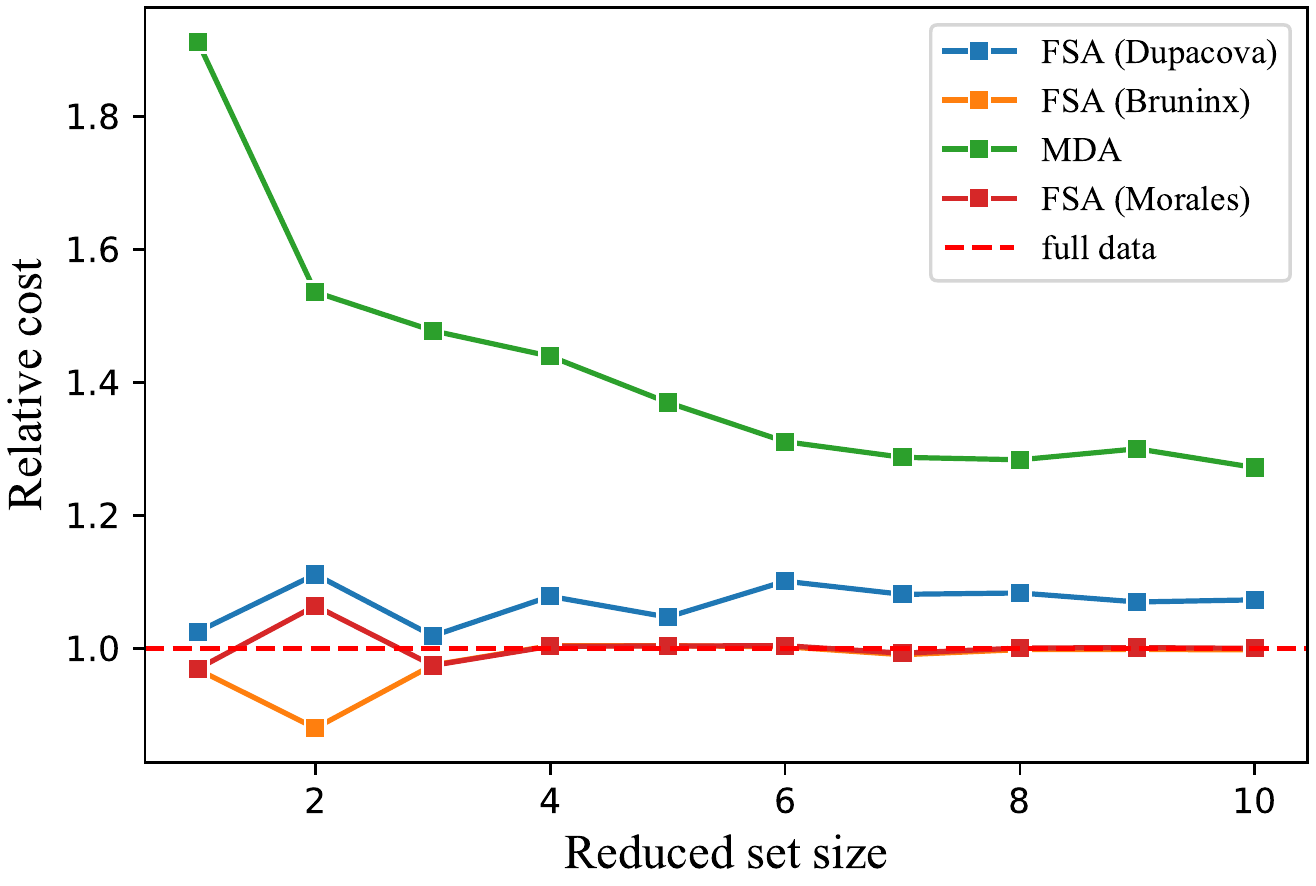}
        \caption{Scenario representation}
    \end{subfigure}
        \hspace{0.05cm}
    \begin{subfigure}{0.32\textwidth}
        \includegraphics[width=\linewidth]{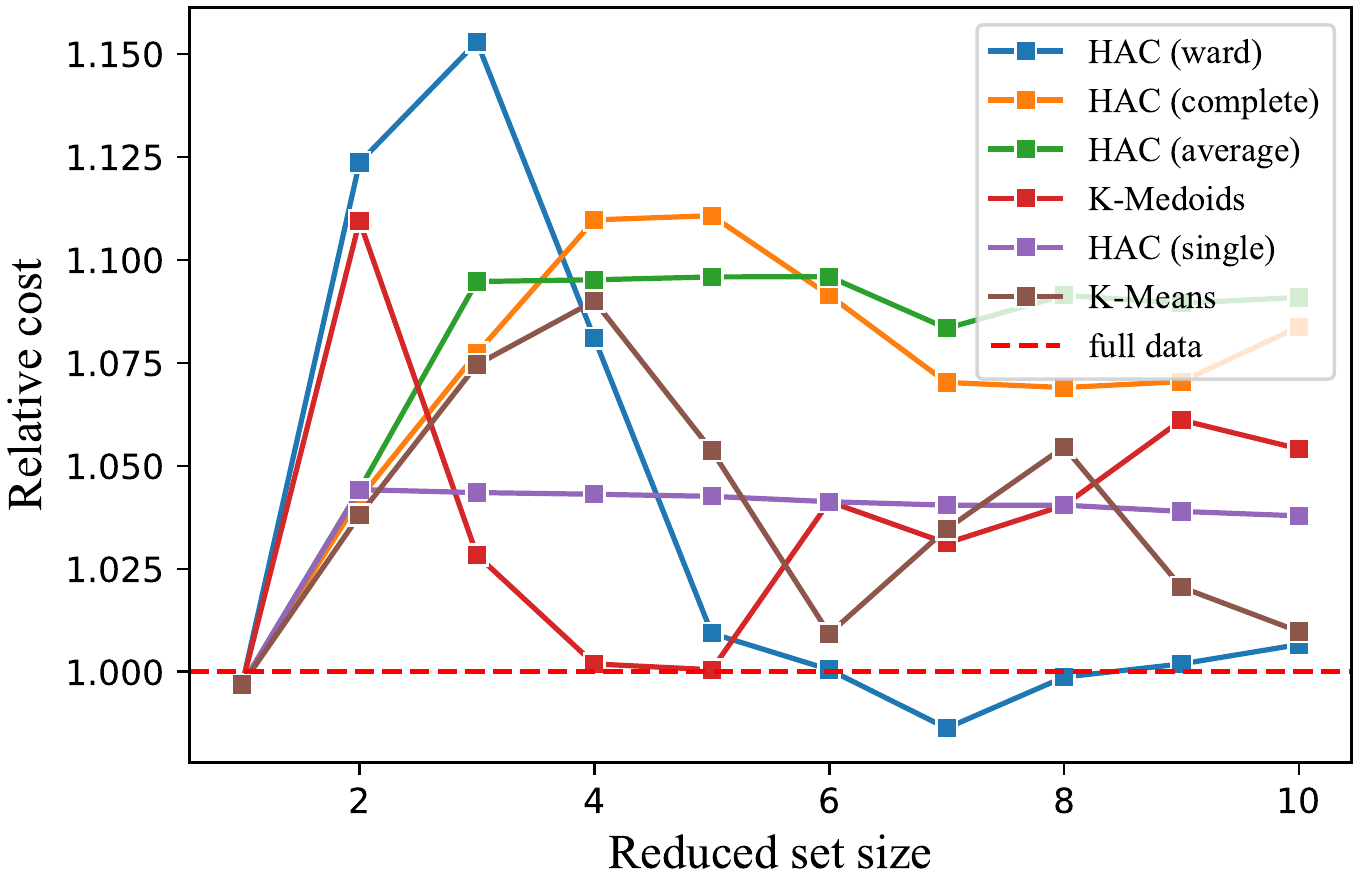}
        \caption{Medoid representation}
    \end{subfigure}
    \caption{In-sample relative objective value change with different reduced set sizes for HRP 38-bus network}
    \label{fig: hrp_insample}
\end{figure}

Just as we did for the Garver system, we plotted the dependence of in-sample relative cost on the reduced set sizes in Figure \ref{fig: hrp_insample}. 
Firstly, the results further confirm some of the explanations given for the Garver system: the centroid representation indeed produces sets that lead to lower bounds for the optimal solution, while solutions obtained with single linkage HAC do not change much with increasing reduced set sizes. Next, our assumption that the jump for MDA in Figure \ref{fig: garver_insample} was connected to data holds true, since that jump is absent for this system. Lastly, the majority of methods with medoid representation are unable to obtain better approximations than the single scenario case, which are equal in all methods with the exceptions of HAC with Ward's linkage and k-medoids.

\begin{table}[h!]
    \centering
        \caption{Comparison of algorithms on HRP 38-bus system}
        \begin{tabular}{ccccc}
            \toprule
            \textbf{Method} & \textbf{OOS gap, \%} & \textbf{Built lines} & \textbf{Capacities, MW} & \textbf{Num. repr.} \\
            \midrule
            All scenarios & - & \{14-27, 15-35, 36-35\} & \{1310.7, 2500, 2500\} & 365 \\
            \midrule
            \multicolumn{5}{c}{Centroid representation} \\
            \midrule
            K-Means & 0.146 & \{14-27, 36-35\} & \{1120, 2500\} & 10 \\
            K-Medoids & 0.143 & \{14-27, 36-35\} & \{1220, 2500\} & 10 \\
            HAC (Ward) & 0.142 & \{14-27, 36-35\} & \{1275.3, 2500\} & 10 \\
            HAC (average) & 0.142 & \{14-27, 36-35\} & \{1350, 2500\} & 10 \\
            HAC (single) & 0.382 & \{14-27, 36-35\} & \{1224.8, 1474.3\} & 10 \\
            HAC (complete) & 0.142 & \{14-27, 36-35\} & \{1295.4, 2500\} & 10 \\
            \midrule
            \multicolumn{5}{c}{Medoid representation} \\
            \midrule
            K-Means & 0.324 & \{14-27, 15-35, 36-35\} & \{646.1, 2500, 1211.4\} & 1 \\
            K-Medoids & 0.070 & \{14-27, 15-35, 36-35\} & \{500, 2500, 2500\} & 5 \\
            HAC (Ward) & 0.070 & \{14-27, 15-35, 36-35\} & \{500, 2500, 2500\} & 6 \\
            HAC (average) & 0.324 & \{14-27, 15-35, 36-35\} & \{646.1, 2500, 1211.4\} & 1 \\
            HAC (single) & 0.324 & \{14-27, 15-35, 36-35\} & \{646.1, 2500, 1211.4\} & 1 \\
            HAC (complete) & 0.324 & \{14-27, 15-35, 36-35\} & \{646.1, 2500, 1211.4\} & 1 \\
            \midrule
            \multicolumn{5}{c}{Scenario representation} \\
            \midrule
            MDA & 0.024 & \{14-27, 15-35, 36-35\} & \{1835.1, 2500, 2500\} & 10 \\
            FSA (Dupačová) & \textbf{$\sim$ 0} & \{14-27, 15-35, 36-35\} & \{1325.2, 2500, 2500\} & 3 \\
            FSA (Morales) & 0.022 & \{14-27, 15-35, 36-35\} & \{1815.7, 2500, 2500\} & 8 \\
            FSA (Bruninx) & 0.035 & \{14-27, 15-35, 36-35\} & \{712.1, 2500, 2500\} & 9 \\
            \bottomrule
        \end{tabular}
    \label{table: hrp_methods_comparison}
\end{table}

\begin{table}[htb]
    \centering
        \caption{Running time comparison of algorithms on HRP 38-bus system}
    \begin{tabular}{ccc}
        \toprule
        \textbf{Method} & \textbf{Method time, s} & \textbf{Total time, s} \\
        \midrule
        All scenarios & - & 14650 \\
        \midrule
        K-Means & $3.67 \pm 0.08$ & 37.7 \\
        K-Medoids & $0.025 \pm 0.003$ & 34 \\
        HAC & $0.297 \pm 0.004$ & 34.3 \\
        MDA & $0.301 \pm 0.002$ & 34.3 \\
        FSA (Dupačová) & $24.1 \pm 0.04$ & 58.1 \\
        FSA (Morales) & $179 \pm 0.15$ & 213 \\
        FSA (Bruninx) & $3144 \pm 12.7$ & 3178 \\
        \bottomrule
    \end{tabular}
    \label{table: hrp_runtime}
\end{table}

Since our methodology for selecting the optimal size of the reduced set relies on the in-sample results to bring it closer to real-life cases, for the majority of medoid methods in Table \ref{table: hrp_methods_comparison} the number of representative scenarios was set to 1, indicating their inability to improve from the starting solution. 
The resulting solution was still able to identify the correct investment decisions but yielded large errors in estimating the required line capacities. 
K-medoids and HAC with Ward's linkage (medoid representation) were able to correctly identify the lines and approximate the capacity for line 36-35 but still yielded a large error for the capacity of line 14-27.
On the other hand, solutions with sets obtained with centroid representation methods were unable to identify all the required investment decisions and omitted line 15-35, but their approximation of line capacity for lines 14-27 and 36-35 was close to the optimal solution with the exception of the HAC with single linkage, which yielded a large error for line 36-35.
Finally, the methods with scenario representation demonstrated the best performance - solutions with these methods were able to identify correct investment decisions and approximate the required line capacities better than medoid methods, whereas the FSA with \ref{eq: Dupacova} cost function produced a solution that had a negligible error when compared to the optimal solution. For the HRP system, a gap of 1\% is equal to roughly \$130 million and investments account for around 1\% of the total cost.

Finally, for all methods we measured TEP optimization time on the HRP system with 10 scenarios. The results for each method can be seen in Table \ref{table: hrp_runtime}. HACs with different linkage criteria again showed similar performance in terms of running time, hence they were combined into a single method. 
The results generally validate the ones obtained for the Garver system: the k-means, k-medoids, HAC, and MDA have roughly equal times, while FSA algorithms are slower. The original FSA is only marginally slower than the other algorithms, while contextual FSAs are considerably slower, especially the FSA with cost function \eqref{eq: Bruninx}, which takes about 50 minutes to compute. 
Considering that contextual cost functions \eqref{eq: Morales} and \eqref{eq: Bruninx} for FSA showed worse performance than the original cost function \eqref{eq: Dupacova}, it makes them less applicable in large-scale systems.

\subsection{Final remarks}

Firstly, the results of our numerical experiments show that the choice of TSA can lead to significant changes in the optimization problem solution. The difference usually lies in the approximation of the line capacities for the built lines, but in some cases, such as with centroid methods in the HRP network, it can lead to incorrect investment decisions in line building. 
Secondly, the clustering algorithms with medoid representation seem to be less robust when there are changes in the network topology, since the majority of them were unable to provide a better solution than the single scenario case for the 38-bus system. 
On the other hand, centroid methods worked for both systems, but in the case of the HRP network, the solutions using those sets were unable to make correct line-building decisions. 
Finally, the FSA with cost function \eqref{eq: Dupacova} and the MDA showed stable and high-quality performance across both tests and had a reasonable computational time. Two other FSA variations also produced good solutions but were significantly slower, making them less preferable to the vanilla FSA.

The sensitivity analysis showed that varying even a single parameter can lead to an almost 8-fold difference in the required computational time. This means that, depending on data, even the small-scale cases might become too expensive to solve, and the problem would become even more significant in large-scale systems. Such behavior in NP-hard problems further stresses the importance of TSA as one of the few reliable avenues towards reducing computational cost in stochastic MILP models.

\section{Conclusion}
\label{sec:Conclusion}

Power systems are facing new challenges in including stochastic generation from the increasing generation shares of renewable sources.  Optimization models for operating and planning power systems are becoming more complex.  During the past decade, hundreds of new aggregation methods for reducing scenario set sizes have been proposed. Concepts from different fields have been borrowed and adapted by the power systems community, such as data clustering or scenario reduction techniques. 

This paper has reviewed the up-to-date developments on scenario aggregation methods for power system optimization problems. Because it is intended to serve as a guide for power system modelers, it describes, classifies and quantitatively compares existing methodologies that comprise the topic's state-of-the-art. 
A general framework for aggregation methods was introduced in this paper. The most common aggregation methods were carefully discussed and classified. A bibliometric analysis highlighting the relevance of the examined works in different fields was presented. Numerical tests were run for the transmission expansion planning problem and discussed.   
From the numerical tests run for time-series scenario aggregation on 16 methods, we conclude that centroid methods could lead to erroneous investment decisions and medoid methods provide a less robust solution for various dimensions of representative scenario sets. On the other hand, methods based on the forward selection algorithm provide more accurate results but there is higher computational overhead.  

As a final remark, we have shown that it is difficult to prescribe a priori the best aggregation method for a wide range of applications. Providing such recommendations for even the same application problem and case study is unachievable when data input varies, which we observed through numerical experiments. However, modelers can find adequate, though not necessarily optimal, aggregation methods that satisfy modelers' requirements on application.

\appendix

\begin{landscape}

\begin{small}
\begin{longtable}{M{0.05\linewidth}M{0.15\linewidth}M{0.3\linewidth}M{0.2\linewidth}M{0.3\linewidth}}
\caption{Summary of aggregation methods}\label{table:Methods}\\
\toprule
\textbf{Reference} & \textbf{Problems} & \textbf{Methods} & \textbf{Objects} & \textbf{Metrics}  \\
\midrule
\endfirsthead

\caption* {Table \ref{table:Methods}: Continued}\\
\toprule
\textbf{Reference} & \textbf{Problems} & \textbf{Methods} & \textbf{Objects} & \textbf{Metrics}  \\
\midrule
\endhead

\endfoot

\endlastfoot

\cite{Palmintier2014} & UC & Heuristic clustering by: unit type, unit type and additional characteristics, plant & Generation units & Operational cost, system-wide CO2 emissions, energy mix, commitment plan, hourly power output \\
\midrule
\cite{Wogrin2014} & UC & K-means with hourly data, 6 load levels or 6 system states & Operational states & Total system cost, power production, electricity prices, reserve prices \\
\midrule
\cite{Poncelet2016} & Long-term planning TIMES; UC LUSYM & Averaged time slices; representative day time slices & Generation and demand scenarios & Approximation RMSE of residual load duration curve (RLDC) for different wind penetration levels \\
\midrule
\cite{Shayesteh2016} & AC and DC OPF, stochastic UC & FSA with \ref{eq: Dupacova} cost function; ATC aggregation & System nodes, generation and demand scenarios & Total cost, total losses, expected energy not served, generator commitment plans \\
\midrule
\cite{TejadaArango2018} & UC & K-medoids; k-means & Generation and demand scenarios & Total cost, generation mix, investment decisions\\
\midrule
\cite{Teichgraeber2019} & Two UC problems: electricity storage and gas turbine dispatch & K-means with centroid and medoid representation; k-medoids; DBA; k-shape; HAC with Ward linkage, Euclidean distance and centroid or medoid representation & Generation and demand scenarios & Total cost \\
\midrule
\cite{Merrick2016} & GEP & Monthly median and peak electricity demand days with 4 hour resolution (SWITCH); mean of all data & Generation and demand scenarios & Duration curves, capacity and generation by technology \\
\midrule
\cite{Liu2018} & GEP & HAC with minmax linkage and DTW distance (HC); k-means for full set, then HC for clusters (kMHC); k-means & Generation and demand scenarios & Peak loads, RMSD for demand, wind and solar generation, RMSD for investment decisions \\
\midrule
\cite{Mallapragada2018} & Two GEPs, UC with ED & Time slices - 4 seasons with 4 times of day; k-means with closest to centroid representation & Generation and demand scenarios & Annual curtailment, unmet demand, annual RES penetration and thermal generation mix \\
\midrule
\cite{Buchholz2019} & GEP & Dummy selection; statistical representation; RLDC selection; dynamic blocking; k-means; kMHC; fuzzy clustering plus HAC with correlation linkage; exhaustive search; optimized RLDC approximation (MILP model); optimized criteria selection (MILP model) & Generation and demand scenarios & Performance index - sum of differences for maximum capacities, generation mixes, invested units, numbers of shut-downs, and total costs \\
\midrule
\cite{Yeganefar2020} & GEP & SOM-based selection of extreme operating points & Generation and demand scenarios & Investment decisions, required dispatchable capacity, operating costs \\
\midrule
\cite{Helist2020} & GEP & Regular decomposition; k-means with closest to centroid representation; random sampling; brute force optimization & Generation and demand scenarios & Total cost, RMSD between duration curves \\
\midrule
\cite{Domnguez2021} & GEP & HAC with Ward linkage (daily and weekly); CTPC; multi-chronological CTPC & Generation and demand scenarios & Total cost, duration curves, total generation capacity, share of unserved demand \\
\midrule
\cite{Gonzato2021} & GEP & MILP optimization methods: ORDF, RR, ORDO; HAC with Ward linkage and medoid representation; CTPC & Generation and demand scenarios & L1-error in approximation of the full solution, mean time series error in approximating time-series values, capacity mix error \\
\midrule
\cite{Fitiwi2015} & TEP & K-means with medoid representation and OPF-based feature engineering & Generation and demand scenarios & Investment costs, total dispatch cost \\
\midrule
\cite{Alvarez2017} & TEP & HAC with average linkage, centroid representation and maximum weighted distance & Generation and demand scenarios & Representativeness of cluster centrums, representativeness of the corridors loading levels \\
\midrule
\cite{Ploussard2017} & TEP & K-means with closest to centroid representation & Operational states & Total system costs, investment decisions \\
\midrule
\cite{MajidiQadikolai2018} & TEP & K-means with scenario bundling (IP optimization model) & Generation and demand scenarios & Optimality gap, total cost, investment decisions \\
\midrule
\cite{Palmintier2011} & Joint GEP and UC & Grouping by similar characteristics & Generation units & Operational cost, energy mix, CO2 emissions \\
\midrule
\cite{Scott2019} & Joint GEP and UC & K-medoids with average squared Euclidean distance & Generation and demand scenarios & Duration curves, error in GEP estimation of costs, additional cost over best found solution, selected technologies and capacities, annual total net demand error \\
\midrule
\cite{Tso2020} & Joint GEP and UC & Iterative decomposition algorithm based on HAC with Ward linkage & Generation and demand scenarios & Total cost \\
\midrule
\cite{Sun2019} & Joint TEP and GEP & HAC with Ward linkage, medoid representation and Euclidean distance & Generation and demand scenarios & Normalized RMSE of total cost \\
\midrule
\cite{vanderHeijde2019} & Multi-node district heating model (modesto) & MILP model for selection of representative days and MIQP model for chronology restoration & Generation portfolios & Energy import, electricity cost, solar energy usage, storage energy loss, pumping energy loss \\
\midrule
\cite{Rsnen2009} & Time-series dimensionality reduction & K-means with feature engineering & Customer load data & Index-of-Agreement \\
\midrule
\cite{SanchezGarcia2014} & Power grid reduction & Spectral clustering with HAC & System nodes & Expansion \\
\midrule
\cite{Vankov2020} & Power grid reduction & Consensus clustering (Spectral + HAC) with Euclidean, MJC, DTW distances & System nodes, generation and demand scenarios & Silhouette score \\
\midrule
\cite{SatreMeloy2020} & Cluster analysis & K-means; HAC with single, complete, average, Ward and centroid linkage criteria & Customer demand data & Calinski-Harabasz score, C-index, Davies-Bouldin index, Duda index, gap statistic, silhouette score \\
\midrule
\cite{Li2020} & Dimensionality reduction & K-means; HAC with complete linkage and Euclidean distance; GMM with EM algorithm; spectral clustering; DBSCAN & Decarbonization pathways & Sum of squared errors, Davies-Bouldin index, Calinski-Harabasz score, Dunn index, silhouette score\\
\midrule
\cite{Tanoto2020} & Cluster analysis & K-means; SOM & Generation portfolios & Calinski-Harabasz score, Davies-Bouldin index, energy trilemma \\
\midrule
\cite{Zhang2021} & Real-time pricing and power consumption scheduling & Iterative decision-making oriented clustering & Demand & Optimality loss, required number of clusters \\
\bottomrule
\end{longtable}
\end{small}

\end{landscape}

 \bibliographystyle{elsarticle-num} 
 \bibliography{cas-refs}





\end{document}